\documentclass{amsart}
\usepackage[mathscr]{eucal}
\usepackage{amssymb, amsmath,array, amscd, enumerate,verbatim,epsfig}
\usepackage[active]{srcltx}
\usepackage{url}

\theoremstyle{change}
\newtheorem{thm}{Theorem}[section]
\newtheorem{THM}{Theorem}
\newtheorem{prop}[thm]{Proposition}
\newtheorem{lemma}[thm]{Lemma}
\theoremstyle{definition}
\newtheorem{definition}[thm]{Definition}

\newtheorem{remark}[thm]{Remark}
\newtheorem{cor}[thm]{Corollary}
\newtheorem{COR}{Corollary}
\newtheorem{problem}[thm]{Problem}

\newcommand{\W}{\ensuremath{\mathcal{W}}}
\newcommand{\F}{\ensuremath{\mathcal{F}}}
\newcommand{\G}{\ensuremath{\mathcal{G}}}

\newcommand{\C}{\ensuremath{\mathbb{C}}}

\newcommand{\oo}{\ensuremath{\mathcal{O}}}
\newcommand{\pp}{\ensuremath{\mathbb{P}}}

\newcommand{\leg}{\operatorname{Leg}}

\def\px{\frac{\partial}{\partial x}}
\def\py{\frac{\partial}{\partial y}}

\title{Rigid flat webs on the projective plane}
\author{David Mar\'{\i}n}
\address{Departament de Matem\`{a}tiques \\ Universitat Aut\`{o}noma de Barcelona \\ E-08193  Bellaterra (Barcelona)\\ Spain \\}
\email{ davidmp@mat.uab.es}

\author{Jorge Vit\'{o}rio Pereira}
\address{Instituto de Matem\'atica Pura e Aplicada\\ Est.
D. Castorina, 110\\
22460-320, Rio de Janeiro, RJ, Brasil\\} \email{ jvp@impa.br}

\date{\today}

\thanks{The first author was partially supported by FEDER / Ministerio de Educaci\'{o}n y Ciencia of Spain, grant MTM 2008-02294. He specially thanks the invitation of IMPA at Rio de Janeiro in August 2009. The second author was partially supported by FAPERJ and Cnpq. He thanks the invitation of CRM at Bellaterra in July 2010.}
\begin{document}

\maketitle

\begin{abstract}
This paper studies global webs on the projective plane with vanishing curvature. The study is based on an interplay of local and
global arguments. The main local ingredient is a criterium for the regularity of the curvature at the neighborhood of a generic
point of the discriminant. The main global ingredient, the Legendre transform, is an avatar of classical projective duality in the
realm of differential equations. We show that the Legendre transform of what we call reduced convex foliations are webs
with zero curvature, and we exhibit  a countable infinity  family of convex foliations which give rise to a family of webs with zero curvature
not admitting non-trivial deformations with zero curvature.
\end{abstract}

\section{Introduction}
Roughly speaking, web geometry is the study of invariants for finite families of foliations. The subject was
initiated by Blaschke and his school in the late 1920's, but among its most
emblematic results there are versions of Lie-Poincar\'{e}-Darboux's converse to Abel's addition
theorem which can be traced back to the XIXth century. While the subject can be developed
in different categories, the earlier practitioners of the subject dealt with finite
families of \emph{germs} of \emph{holomorphic} foliations.

Recently, the study of holomorphic webs globally defined on compact complex manifolds
started to be pursued, see for instance \cite{Yartey, CL1,  MPP, PP}.  It is in this
context that this work places itself. Its main purpose is to investigate the irreducible
components of the space of flat webs on the projective plane.

\subsection{Webs on the projective plane}
In the same way that a foliation  on the projective plane is defined by a polynomial
$1$-form  $a(x,y) dx + b(x,y) dy$ on $\mathbb C^2$ with isolated zeros, a $k$-web
on the projective plane is defined by a $k$-symmetric polynomial $1$-form
\[
\omega = \sum_{i+j = k} a_{ij}(x,y) dx^i dy^j
\]
with isolated zeros and non identically zero discriminant. In more intrinsic
terms, a $k$-web on a complex surface $S$ is defined by an element $\omega$ of
$H^0(S,Sym^k \Omega^1_S\otimes N)$ for a suitable line-bundle $N$, still
subjected to the two conditions above: isolated zeros and non-zero discriminant.

When $S = \mathbb P^2$, it is natural to write $N$ as $\mathcal O_{\mathbb P^2}(d + 2k)$
since the pull-back of $\omega$ to a  line $\ell \subset \mathbb P^2$ will
be a section of $Sym^k \Omega^1_{\mathbb P^1} ( d+ 2k) = \mathcal O_{\mathbb P^1}(d)$
and consequently for a generic $\ell$ the integer  $d$  will count the number of
tangencies between $\ell$ and the $k$-web $\mathcal W$ defined by $\omega$. That said, we promptly
see that $\mathbb W(k,d)$ -- the space of $k$-webs on $\mathbb P^2$ of degree $d$ -- is an open subset of
$\mathbb P H^0(\mathbb P^2, Sym^k \Omega^1_{\mathbb P^2}(d+2k))$.

\subsection{Curvature and flatness}\label{curvflat} One of the first results of web geometry, due to Blaschke-Dubourdieu,  characterizes the
local equivalence of a (germ of)  $3$-web $\mathcal W$ on $\mathbb C^2$ with the {\it trivial} $3$-web defined
by $dx \cdot dy \cdot (dx-dy)$ through the vanishing of a differential covariant: the curvature of $\mathcal W$.
It is a meromorphic $2$-form $K(\mathcal W)$ with poles on the discriminant of $\mathcal W$ that satisfies
$\varphi^* K ( \mathcal W) = K ( \varphi^* \mathcal W)$ for any bihomolomorphism $\varphi$.

For a $k$-web $\mathcal W$ with $k>3$, one usually defines the curvature of $\mathcal W$ as the
sum of the curvatures of  all $3$-subwebs of $\mathcal W$. It is again a differential covariant, and to the
best of our knowledge there is no result characterizing its  vanishing (a conjectural characterization for $4$-webs is proposed in \cite{RipollT}).  Nevertheless,  according to a result of Mihaileanu -- recently rediscovered by H\'{e}naut, Robert, and Ripoll
-- this vanishing is a necessary condition for the maximality of the rank of the web, see \cite{Henaut,Ripoll,bourbaki}  for a
thorough discussion and   pertinent references.

\medskip

The $k$-webs with zero curvature are here called {\bf flat $k$-webs}, and the subset of $\mathbb W(k,d)$ formed
by the flat $k$-webs will be denoted by $\mathbb{FW}(k,d)$. It is  a Zariski closed subset of $\mathbb W(k,d)$
and it is our purpose to describe some of its  irreducible components. More specifically we will characterize one
irreducible component of $\mathbb{FW}(k,1)$ for each $k \ge 3$. For that sake we will pursue the following strategy:
\begin{enumerate}
\item study the {\bf regularity of the curvature} on irreducible components of the discriminant;
\item translate constraints imposed by (1) on  flat $k$-webs of degree $1$ into
constraints on foliations of degree $k$ using {\bf projective duality};
\item apply (2) to  {\bf convex foliations} to establish the   flatness of their duals;
\item apply (1) combined with (2) to determine the deformations of convex foliations with flat duals.
\end{enumerate}

We will now proceed to a more detailed discussion about each of the steps of our strategy, and will take
the opportunity to state the  main results of this work.

\subsection{Regularity of the curvature} As mentioned above,
the curvature of a web $\mathcal W$ on a complex surface is a meromorphic $2$-form with poles contained in the discriminant $\Delta(\mathcal W)$ of
$\mathcal W$.
As there are no  holomorphic $2$-forms on
$\mathbb P^2$,
 the curvature of a global  web  $\mathcal W $ on the projective plane is zero if and only if it is holomorphic over the generic points of the irreducible components of $\Delta(\mathcal W)$.

This very same observation was used in \cite{CDQL} to classify   completely decomposable quasi-linear (CDQL) exceptional webs  on $\mathbb P^2$ with zero curvature. There, a criterium  for the holomorphicity of the curvature over an irreducible component
of $\Delta(\mathcal W)$ is given under a certain number of hypothesis. Among these hypothesis, there is the local decomposability of $\mathcal W$, that is $\mathcal W$ can be locally written as a product of foliations. While this was sufficient in that setup, here we will deal with webs which are not necessarily locally decomposable.

\medskip

If $\mathcal W$ is a germ of $(k+2)$-web on $(\mathbb C^2,0)$ with reduced,  smooth, and non-empty  discriminant $\Delta(\W)$ then
it is the superposition of an irreducible $2$-web $\mathcal W_2$ and a completely decomposable web $\mathcal W_k$.
Moreover $\Delta(\mathcal W_2) = \Delta(\mathcal W)$ and $\Delta(\mathcal W_k) = \emptyset$. Our
first result is a generalization of  \cite[Theorem 7.1]{CDQL}, with hypothesis also satisfied by webs with
reduced discriminant.

\begin{THM}\label{T:1}
Let $\mathcal W$ be a germ of $(k+2)$-web on $(\mathbb C^2,0)$ with smooth (but not necessarily reduced), and non empty discriminant.
Assume $\mathcal W= \mathcal W_2 \boxtimes \mathcal W_k$  where $\mathcal W_2$ is a $2$-web satisfying $\Delta(\mathcal W_2)= \Delta(\mathcal W)$,
and $\mathcal W_k$ is a $k$-web.
The curvature of $\W$ is holomorphic along $\Delta(\mathcal W)$ if and only if $\Delta(\W)$ is invariant by either $\W_{2}$ or  $\beta_{\W_{2}}(\W_{k})$.
\end{THM}

In the statement $\beta_{\mathcal W_2}(\mathcal W_k)$ stands for the $\mathcal W_2$-barycenter of $\mathcal W_k$. It is a $2$-web naturally associated
to the pair $(\mathcal W_2, \mathcal W_k)$ as defined in Section \ref{S:bary}.

\medskip

A simple consequence of Theorem \ref{T:1} is the following result, which will play an essential role in our
study of irreducible components of $\mathbb F \mathbb W(k,1)$.

\begin{COR}\label{cor:1}
Let $\W=\W_{2}\boxtimes\W_{k}$ be a $(k+2)$-web  in $(\C^{2},0)$ such that $\Delta(\W_{2})=\Delta(\W)$ is smooth and invariant by $\W_{2}$. Then any deformation $\W^{\varepsilon}$ of $\W$ having holomorphic curvature is of the form $\W^{\varepsilon}=\W_{2}^{\varepsilon}\boxtimes\W_{k}^{\varepsilon}$ with $\Delta(\W_{2}^{\varepsilon})=\Delta(\W^{\varepsilon})$ invariant by $\W_{2}^{\varepsilon}$.
\end{COR}

\subsection{Legendre transform} Browsing classical books on ordinary differential equations one
can find the so called Legendre transform, see for instance \cite[page 40]{Ince}. It is an involutive transformation
which sends the polynomial differential equation   $F(x,y,p)$ to $F(P,XP-Y,X)$, where $p= dy/dx$ and $P=dY/dX$.
It can be expressed in global projective coordinates, as Clebsch already did back in the XIXth century \cite{Darboux} and as we explain in Section \ref{S:legendre}.
It turns out to be an isomorphism between $H^0(\mathbb P^2, Sym^k \Omega^1_{\mathbb P^2}(d+ 2k))$ and $H^0(\mathbb P^2, Sym^d \Omega^1_{\mathbb P^2}(k+ 2d))$,
and as such associates to a   $k$-web of degree $d$, a $d$-web of degree $k$.

\smallskip

There is a  beautiful underlying geometry  which we take our time to discuss. We analyze carefully the dual of foliations.
Radial singularities and invariant components of the inflection curve turn out to have a distinguished behavior.
Looking at foliations with extremal properties  with respect to the latter we are  able  to put in evidence an
infinite family of examples of webs with zero curvature.

\subsection{Convex foliations} More precisely, we look at the dual of what we call {\bf reduced  convex foliation}. For
us a foliation $\mathcal F$  on $\mathbb P^2$ is convex if its leaves other than straight lines have no inflection points. In other words,
the inflection divisor $I(\mathcal F)$ of $\mathcal F$ -- called in \cite{extactic} the first extactic divisor --
is completely invariant by $\mathcal F$. When besides being completely invariant, this divisor is also reduced
we will say that $\mathcal F$ is a reduced convex foliation.

Our second main result is about the dual of reduced convex foliations and it can be phrased as follows.

\begin{THM}\label{T:convex}
If $\mathcal F$ is a reduced convex foliation of degree $d\ge 3$ then its Legendre transform is a flat $d$-web of degree one.
\end{THM}

Of course such result would be meaningless if examples of reduced convex foliations did not exist.  Fortunately,
this is far to be true as we have for every $d\ge 2$, the reduced convex foliation $\mathcal F_d$ of degree $d$
defined by the levels of the rational function $\frac{ x^{d-1}(y^{d-1} - z^{d-1})}{  y^{d-1}(x^{d-1} - z^{d-1}) }$.
It turns out that the dual webs are not just flat but indeed algebraizable, see Proposition \ref{P:fermat}.

\medskip

Besides the infinite family $\mathcal F_d$ we are aware of three other examples of reduced convex foliations. The Hesse pencil
of degree four; the Hilbert modular foliation of degree $5$ studied in \cite{Hilb}; and one foliation of degree $7$ induced
by a pencil of curves of degree $72$ and  genus $55$ related to extended Hesse arrangement. These examples are described in
Section \ref{S:examples}. In Table \ref{Table:1} we list the number of radial singularities of these
examples and their main birational invariants: Kodaira and numerical Kodaira dimension as defined in \cite{canonical}, see also \cite{Br1, LG}.

\begin{table}[ht]
\centerline{\begin{tabular}{|l|c|c|c|c|c|l|} \hline
{\bf Fol.} & $d(\F)$ & $r(\F)$ & {\rm kod($\mathcal{F}$)} &  $\nu(\mathcal{F})$ &  {\bf description} \\
\hline
$\F_2$   & $2$ & $4$ & $-\infty$ & $-\infty$ &   \small{rational fibration}\\
\hline
$\F_3$   & $3$ & $7$ & $-\infty$  & $-\infty$ & \small{rational fibration}\\
\hline
$\F_4$   & $4$ & $12$ & $0$ & $0$ &      \small{isotrivial elliptic fibration}\\
\hline
$\F_d$    & $d \ge 5$ & $(d-1)^2 + 3$ & $1$ & $1$ &  \small{isotrivial hyperbolic fibration}\\
\hline
$\mathcal H_{4}$  & $4$ & $9$ & $1$ & $1$ &  \small{non-isotrivial elliptic fibration}\\
\hline
$\mathcal H_5$ & $5$ & $16$ & $-\infty$ & $1$ &  \small{Hilbert Modular foliation} \\
\hline
$\mathcal H_{7}$ & $7$ & $21$ & $2$ & $2$ &   \small{non-isotrivial hyperbolic fibration} \\
\hline
\end{tabular}}
\medskip

\caption{Known examples of reduced convex foliations.}\label{Table:1}
\end{table}

\subsection{Rigidity} Our third main result concerns the deformations of the webs dual to the foliations $\mathcal F_d$ inside $\mathbb F \mathbb W(d,1)$.
It can be succinctly stated as follows.

\begin{THM}\label{T:3}
If  $d= 3$ or $d \ge 5$  then the closure of the $PGL(3,\mathbb C)$-orbit of the Legendre transform of  ${\mathcal F}_d$ is an
irreducible component of $\mathbb F\mathbb W(d,1)$. For $d=4$, the closure of the $PGL(3,\mathbb C)$-orbit
of the Legendre transform of ${\mathcal F}_4$ has codimension one in an irreducible component of $\mathbb F\mathbb W(4,1)$.
\end{THM}

Indeed we prove slightly more, as we describe a Zariski open subset of the irreducible
component of $\mathbb F\mathbb W(4,1)$ containing the Legendre transform of ${\mathcal F}_4$.

\subsection{Acknowledgments} We are grateful to Olivier Ripoll. In an early stage of this project we made an extensive use of Ripoll's Maple scripts
to compute the H\'{e}naut's curvature of $4$ and  $5$-webs. The period of experimentation with Ripoll's script was essential as it helped  to
build our  intuition on the subject.  We are also grateful to Maycol Falla Luza for pointing out a number of misprints and  mistakes in  previous versions of this work.

\section{Regularity of the curvature}

Let $\W$ be a web and let $C\subset\Delta(\W)$ be an irreducible component of its discriminant. We say that $C$ is invariant (resp. totally invariant) by $\W$ if and only if $TC\subset T\W|_{C}$ (resp. $TC=T\W|_{C}$) over the regular part of $C$. Notice that when $\W$ is a germ of irreducible web then the two notions coincide.

\subsection{Barycenters of webs}
\label{S:bary}
Theorem \ref{T:1} was proved in \cite{CDQL} in the case that $\W_{2}$ is reducible, so we only need to show it when $\W_{2}$ is irreducible. To this end, we recall and slightly extend the notion of barycenters of webs introduced there.

Let $\W$ be a $k$-web on a complex surface $S$ and let $\F$ be a foliation transverse to $\W$ at
some open set $U\subset S$. For each point $p\in U$ the tangent lines of $\W$ at $p$ can
be considered as $k$ points in  the affine line $\mathbb{P}T_{p}U\setminus [T_{p}\F]$.
Thus, we can consider their barycenter. As $p$ varies on $U$ we obtain
a line distribution, which determines  a  foliation $\beta_{\F}(\W)$ on $U$,
called the \emph{barycenter} of $\W$ with respect to $\F$.
Taking a suitable  system of coordinates in a neighborhood $U$ of each point we can identify $\F$ and $\W$ with its respective slopes
$f,w_{1},\dots,w_{k}:U\to\C$. If we consider the polynomial $W(x):=\prod\limits_{i=1}^{k}(x-w_{i})$ then $\beta_{\F}(\W)$ corresponds to the slope $f-\frac{k\,W(f)}{W'(f)}:U\to\C\cup\{\infty\}$. We note that when $\F$ and $\W$ are not transverse $W(f)=0$ and $\beta_{\F}(\W)$ has slope $f$, even in the case that $W(f)=W'(f)=0$.
We extend the definition of barycenter by replacing the center foliation $\mathcal F$ by
a center  $\ell$-web $\W'$. The extension is straightforward, if  we write pointwise $\W'=\F_{1}\boxtimes\cdots\boxtimes\F_{\ell}$ then
we define the $\mathcal W'$-barycenter of $\mathcal W$ as being  $\beta_{\W'}(\W)=\beta_{\F_{1}}(\W)\boxtimes\cdots\boxtimes\beta_{\F_{\ell}}(\W)$.

\subsection{Curvature}\label{SS:curvatura}
Let us recall the
definition of curvature for a $k$-web $\W$. Let us first assume that $\W$ is a germ of
completely decomposable $k$-web $
{\mathcal W}=
 \mathcal F_1\boxtimes \cdots \boxtimes
\mathcal F_k$. We start by considering $1$-forms $\omega_i$ with
isolated singularities such that $\mathcal F_i = [\omega_i]$.
Following \cite{CDQL}, for every  triple $(r,s,t)$ with $1\le r<s<t\le k$ we define
\[
\eta_{rst} = \eta(\mathcal F_r \boxtimes \mathcal F_s \boxtimes
\mathcal F_t)
\]
as the unique {meromorphic} $1$-form such that
\[
\displaystyle{\left\lbrace \begin{array}{lcl} d
(\delta_{st}\, \omega_r) &=&
\eta_{rst} \wedge \delta_{st}\, \omega_r \\
d (\delta_{tr}\, \omega_s) &=&
\eta_{rst} \wedge \delta_{tr}\, \omega_s  \\
d (\delta_{rs}\,\omega_t) &=&
\eta_{rst} \wedge \delta_{rs}\, \omega_t  \\
\end{array} \right. }
\]
where the function $\delta_{ij}$  is characterized by the relation \[\omega_i \wedge \omega_j =
\delta_{ij}\,dx\wedge dy .\]

Although the $1$-forms $\omega_i$ are not uniquely defined, the $1$-forms
$\eta_{rst}$ are well-defined modulo the addition of a closed
holomorphic $1$-form. The {\it curvature} of the web $\mathcal W= \mathcal F_1\boxtimes \cdots \boxtimes \mathcal F_k$ is
defined by the formula
\[
K({\mathcal W})=
K(\mathcal F_1\boxtimes \cdots \boxtimes
\mathcal F_k) = d \, \eta({\mathcal W})
\]
where
$
\eta({\mathcal W})=
\eta(\mathcal F_1\boxtimes \cdots \boxtimes
\mathcal F_k) = \sum_{1\le r<s<t\le k}  \eta_{rst} \, .
$
It can be checked that $K({\mathcal W})$ is a {meromorphic} $2$-form intrinsically attached to ${\mathcal W}$. More precisely
 for any dominant holomorphic map $\varphi$, one has
$K(\varphi^*{\mathcal W})=\varphi^* \big(K({\mathcal W})  \big).$
This property of the  curvature allows us to extend the definition of curvature to an arbitrary (not necessarily completely decomposable) $k$-web.
If we pass to a ramified Galois covering where the web becomes completely decomposable then the curvature of this new web turns out to
be invariant by the action of the Galois group and descends to a meromorphic $2$-form on our original surface.

\subsection{Lemmata} We now establish some preliminary results aiming at the proofs of Theorem \ref{T:1} and Corollary \ref{cor:1}.
The first is a normal form for germs of $2$-webs with smooth discriminant which is not invariant.

\begin{lemma}\label{nf2}
Let $\W_{2}$ be an irreducible $2$-web and let $C\subset\Delta(\W_{2})$ be a smooth  irreducible component non invariant by $\W_2$. Then there exist a local coordinate system $(U,(x,y))$ such that $C\cap U=\{y=0\}$ and
$\W_{2}|_{U}$ is given by $dx^{2}+y^{m}dy^{2}$, for some odd positive integer $m$.
\end{lemma}
\begin{proof} We mimic the proof given in \cite[\S 1.4]{Arnold} when $m=1$.
First, we can write locally $C=\{y=0\}$ and $\W_{2}: dx^{2}+y^{m}(2\alpha(x,y)dx\,dy+\beta(x,y)dy^{2})=0$ by redressing
the distributions of lines $T\W_{2}|_{C}$ along $C$. Since the foliation $\beta_{dy}(\W_{2}):dx+y^{m}\alpha(x,y)dy=0$ is non singular there exists a function $\tilde x(x,y)$ transverse to $y$ such that $\beta_{dy}(\W_{2})$ is defined by the differential form $d\tilde x$. Writing $\W_{2}:d\tilde x^{2}+y^{m}(2\tilde\alpha(\tilde x,y)d\tilde x\,dy+\tilde{\beta}(\tilde x,y)dy^{2})=0$ and $\beta_{dy}(\W_{2}):d\tilde x+y^{m}\tilde\alpha(\tilde x,y)dy=0$ in the coordinates $(\tilde x,y)$, we deduce that $\tilde\alpha=0$. Thanks to the irreducibility of $\W_{2}$ we can assume that it is given by $dx^{2}+y^{m}\beta(x,y)dy^{2}=0$, where $m$ is odd and $y\not|\beta$.
Taking the pull-back by the ramified covering $\bar y\mapsto y=\bar y^{2}$, there is a unity $u$ such that
$$dx^{2}+4\bar y^{2(m+1)}\beta(x,\bar y^{2})d\bar y^{2}=(dx+\bar y^{m+1}u(x,\bar y^{2})d\bar y)(dx-\bar y^{m+1}u(x,\bar y^{2})d\bar y).$$
There also exists a unity function $v(x,\bar y)$ such that $d(x\pm\bar y^{m+2}v(x,\bar y))$ is parallel to $dx\pm\bar y^{m+1}u(x,\bar y^{2})d\bar y$. Write $v(x,\bar y)=w(x,\bar y^{2})+\bar y z(x,\bar y^{2})$ and define
$$\hat x:=x+\bar y^{m+3}z(x,\bar y^{2})\quad\textrm{and}\quad\tilde y:=\bar yw(x,\bar y^{2})^{\frac{1}{m+2}},$$ by using that $w$ is a unity. Finally return downstairs by putting $$\hat y:=\tilde y^{2}=y w(x,y)^{\frac{2}{m+2}}$$ and verifying that $$\hat x=x+y^{\frac{m+3}{2}}z(x,y)$$ is a well defined change of coordinates because $m$ is odd.
Since
$$x\pm\bar y^{m+2}v(x,\bar y)=\hat x\pm\tilde y^{m+2}=\hat x\pm\hat y^{\frac{m+2}{2}},$$
we deduce that
$\W_{2}:d\hat x^{2}-\left(\frac{m+2}{2}\right)^{2}\hat y^{m}d\hat y^{2}=0$ which can be reduced to the normal form $dx^{2}+y^{m}dy^{2}$ by rescaling.
\end{proof}
The second preliminary result provides an asymptotic expansion of the curvature of a decomposable $3$-web along an irreducible component (which we can assume to be $y=0$) of its discriminant.
\begin{lemma}\label{curv}
Let $0\le a_{1}\le a_{2}\le a_{3}$ be integers and consider  the $3$-web $\W$ defined by the  $1$-forms
\[
 \omega_{i}=dx+y^{a_{i}}h_{i}(x,y)dy \, , \quad i=1,2,3 \, .
\]
Suppose that the functions $h_1, h_2,$ and $h_3$ do not vanish along $\{y=0\}$, and that the same holds true
for the differences $h_{i}-h_{j}$ when  $a_{i}=a_{j}$ for $i\neq j$. Then
the  curvature of $\W$ is the exterior differential of the meromorphic $1$-form
$$\left(\frac{a_{1}-a_{2}}{h_{31}(x,0)}\frac{1}{y^{a_{1}+1}}+\left[\frac{h_{23}\partial_{y}h_{12}-h_{12}\partial_{y}h_{23}}{h_{12}h_{23}h_{31}}\Big|_{y=0}\right]\frac{1}{y^{a_{1}}}+\cdots\right)dx+\left(\frac{a_{1}}{y}+\cdots\right)\,dy$$
where the dots correspond to  higher order terms in the variable $y$ and
$$h_{ij}=\left\{\begin{array}{lcl} \hphantom{h_{j}}-h_{i}  & \text{if} & a_{i}<a_{j}, \\
h_{j}-h_{i} & \text{if} & a_{i}=a_{j},\\
h_{j}  & \text{if} & a_{i}>a_{j}.
 \end{array}\right.$$
\end{lemma}

\begin{proof}
We will use the notations of Section~\ref{SS:curvatura}. Notice that under our assumptions
$\delta_{ij}=y^{a_{i}}h_{i}-y^{a_{j}}h_{j}$. If we write $\eta=A\,dx+B\,dy$ then
the equalities
$$d(\delta_{ij}\omega_{k})=\big(-\partial_{y}\delta_{ij}+y^{a_{k}}\partial_{x}(\delta_{ij}h_{k})\big)dx\wedge dy=\delta_{ij}\big(A y^{a_{k}}h_{k}-B\big)dx\wedge dy=\eta\wedge\delta_{ij}\omega_{k},$$
where $(i,j,k)$ runs over the cyclic permutations of $(1,2,3)$, are equivalent to the linear system
$$\left(\begin{array}{cc} \delta_{12}y^{a_{3}}h_{3} & -\delta_{12}\\
\delta_{23}y^{a_{1}}h_{1} & -\delta_{23}\end{array}\right)\left(\begin{array}{c}A\\ B\end{array}\right)=\left(\begin{array}{c}
\partial_{x}(\delta_{12}h_{3})y^{a_{3}}-\partial_{y}\delta_{12}\\
\partial_{x}(\delta_{23}h_{1})y^{a_{1}}-\partial_{y}\delta_{23}\end{array}\right).$$
The determinant of the system is $\delta=\delta_{12}\delta_{23}(y^{a_{1}}h_{1}-y^{a_{3}}h_{3})=\delta_{12}\delta_{23}\delta_{31}=y^{2a_{1}+a_{2}}h_{12}h_{23}h_{31}+\cdots$.
Since $h_{12}h_{23}h_{31}$ is not a multiple of $y$, it follows that $\delta$ has order $2a_{1}+a_{2}$ at $y=0$.
Solving the system by Cramer's rule we obtain that
\begin{eqnarray*}\delta A&=&\left|\begin{array}{cc}
\partial_{x}(\delta_{12}h_{3})y^{a_{3}}-\partial_{y}\delta_{12} & -\delta_{12}\\
\partial_{x}(\delta_{23}h_{1})y^{a_{1}}-\partial_{y}\delta_{23} & -\delta_{23}
\end{array}\right|\\ &=&y^{a_{1}+a_{2}-1}(a_{1}-a_{2})h_{23}h_{12}+y^{a_{1}+a_{2}}(h_{23}\partial_{y}h_{12}-h_{12}\partial_{y}h_{23})+\cdots
\end{eqnarray*}
and consequently
$$ A=\frac{a_{1}-a_{2}}{h_{31}}\frac{1}{y^{a_{1}+1}}+\frac{h_{23}\partial_{y}h_{12}-h_{12}\partial_{y}h_{23}}{h_{12}h_{23}h_{31}}\frac{1}{y^{a_{1}}}+\cdots$$
On the other hand,
\begin{eqnarray*}
\delta B&=&\left|\begin{array}{cc}
\delta_{12}y^{a_{3}}h_{3}&\partial_{x}(\delta_{12}h_{3})y^{a_{3}}-\partial_{y}\delta_{12}\\
\delta_{23}y^{a_{1}}h_{1}&\partial_{x}(\delta_{23}h_{1})y^{a_{1}}-\partial_{y}\delta_{23}
\end{array}\right|\\&=&a_{1}y^{2a_{1}+a_{2}-1}h_{1}h_{12}h_{23}-a_{2}y^{a_{1}+a_{2}+a_{3}-1}h_{12}h_{23}h_{3}+\cdots,
\end{eqnarray*}
so that $B-\frac{a_{1}}{y}$ is holomorphic along $y=0$.
\end{proof}

Our next preliminary result settles Theorem \ref{T:1} when the discriminant is not invariant.

\begin{lemma}\label{noninv2}
Let $\W_{2}$ be an irreducible $2$-web and let $C\subset\Delta(\W_{2})$ be a smooth  irreducible component non invariant by $\W_2$. Let $\W_{d-2}$ be a smooth web transverse to $\W_2$ along $C$. Then the curvature of $\W_{2}\boxtimes\W_{d-2}$ is holomorphic along $C$ if and only if $C$ is invariant by the barycenter $\beta_{\W_{2}}(\W_{d-2})$.
\end{lemma}
\begin{proof}  We use the normal form  for $\W_{2}$ given by Lemma \ref{nf2} and we write $\W_{d-2}$ as $\prod\limits_{i=1}^{d-2}(dy+c_{i}(x,y)dx)=0$. After passing to the double cover $\pi(x,y)=(x,y^{2})$ we can write
$\pi^*(\W_{2}\boxtimes\W_{d-2})$ as
$$(dx-y^{m+1}dy)(dx+y^{m+1}dy)\prod\limits_{i=1}^{d-2}(c_{i}(x,y^{2})dx+2ydy)=0.$$
Its curvature is the exterior differential of
$$\sum_{1\le i\le d-2}\eta_{+,-,i}+\sum_{\varepsilon=\pm}\sum_{1\le i<j\le d-2}\eta_{\varepsilon,i,j}+\sum_{1\le i<j<k\le d-2}\eta_{ijk}.$$
The last summatory is holomorphic along $y=0$ because it is the pull-back by $\pi$ of the $1$-form $\eta(\W_{d-2})$ associated to the smooth web $\W_{d-2}$.
Writing $\omega_{\pm}=dx\pm y^{m+1}dy$ and $\omega_{i}=dx+\frac{2y}{c_{i}(x,y^{2})}dy$, $i=1,\ldots,d-2$,
we can apply Lemma~\ref{curv} to deduce that $\eta_{+,-,i}+\frac{m}{2}c_{i}(x,0)\frac{dx}{y^{2}}-\frac{dy}{y}$ and $\eta_{\pm,i,j}-\frac{dy}{y}$ are also holomorphic along $y=0$. We conclude
that the curvature of $\pi^*(\W_{2}\boxtimes\W_{d-2})$ is the exterior differential of
$$-\frac{m}{2}\left(\sum_{i=1}^{d-2}c_{i}(x,0)\right)\frac{dx}{y^{2}}+(d-2)^{2}\,\frac{dy}{y}+\eta,$$
where
$\eta$ is a holomorphic $1$-form along $y=0$. Consequently, the curvature of $\W_{2}\boxtimes\W_{d-2}$ is holomorphic along $y=0$ if and only if $\sum_{i=1}^{d-2}c_{i}(x,0)=0$. On the other hand, the restriction of the barycenter $\beta_{\W_{2}}(\W_{d-2})$ to $y=0$ is given by
$dy+\frac{1}{d-2}\left(\sum_{i=1}^{d-2}c_{i}(x,0)\right)dx$, and consequently $y=0$ is invariant by it if and only if $\sum_{i=1}^{d-2}c_{i}(x,0)=0$.
\end{proof}

Finally we deal with the case of invariant discriminant.

\begin{lemma}\label{inv2}
Let $\W_{2}$ be an irreducible $2$-web and let $C\subset\Delta(\W_{2})$ be a smooth  irreducible component of its discriminant  invariant by  $\W_{2}$. Let $\W_{d-2}$ be a smooth web transverse to $C$. Then the curvature of $\W_{2}\boxtimes\W_{d-2}$ is holomorphic along $C$.
\end{lemma}
\begin{proof} If $C=\{y=0\}$ then $\W_{2}$ can be presented by $dy^{2}+y^{m}\eta dx$, for some $1$-form $\eta$ and some integer $m\ge 1$. Reasoning as in the beginning of the proof of Lemma~\ref{nf2} we can assume that $\eta$ is proportional to $dx$.
By passing to the double cover $\pi(x,y)=(x,y^{2})$, we obtain that $\pi^*(\W_{2}\boxtimes\W_{d-2})=\F_{-}\boxtimes \F_{+}\boxtimes\F_{1}\boxtimes\cdots\boxtimes\F_{d-2}$, where $\F_{\pm}:dy\pm y^{m-1}fdx=0$ and $\F_{i}|_{y=0}:dx=0$. The curvature of $\pi^*(\W_{2}\boxtimes\W_{d-2})$ is the sum
$$\sum_{i=1}^{d-2}K(\F_{-}\boxtimes\F_{+}\boxtimes\F_{i})+\sum_{i<j}\sum_{\varepsilon=\pm}K(\F_{\varepsilon}\boxtimes\F_{i}\boxtimes\F_{j})+\sum_{i<j<k}K(\F_{i}\boxtimes\F_{j}\boxtimes\F_{k}).$$
The first term is holomorphic thanks to Theorem \ref{T:1} (in the decomposable case already proved in \cite{CDQL}) because $y=0$ is $\F_{\pm}$-invariant.
 The second term is holomorphic also by Theorem \ref{T:1}. To see that, we shall distinguish two cases. If $m=1$ then $\sum\limits_{\varepsilon=\pm}K(\F_{\varepsilon}\boxtimes\F_{i}\boxtimes\F_{j})$ is the curvature of the $4$-web $\F_{+}\boxtimes\F_{-}\boxtimes\F_{i}\boxtimes\F_{j}$ whose discriminant $y=0$ is invariant by $\beta_{\F_{i}}(\F_{+}\boxtimes\F_{-})=\beta_{dx}(dy^{2}-f^{2}dx^{2})=dy$. If $m>1$ then
 $y=0$ is invariant by $\beta_{\F_{i}}(\F_{\varepsilon})=\F_{\varepsilon}$.  Finally the third term is holomorphic because it is equal to $\pi^*K(\W_{d-2})$.
\end{proof}

\subsection{Proofs of Theorem \ref{T:1} and Corollary \ref{cor:1}} Now we have just to put the previous results together to obtain  proofs of Theorem \ref{T:1}
and Corollary \ref{cor:1}.

\begin{proof}[Proof of Theorem \ref{T:1}]
As we have already mentioned, we restrict to the case that $\W_{2}$ is irreducible.
If $C\subset\Delta(\W_{2})$ is invariant by $\W_{2}$ or by $\beta_{\W_{2}}(\W_{d-2})$ then the curvature of $\W_{2}\boxtimes\W_{d-2}$ is holomorphic along $C$ thanks to Lemmas \ref{inv2} and \ref{noninv2}. Suppose now that  $K(\mathcal W_{2} \boxtimes \mathcal W_{d-2})$ is a  holomorphic $2$-form at a
neighborhood of a generic point of  $C$.
If $C$ is not $\mathcal W_2$-invariant then Lemma 2.2 implies $C$ is invariant by the barycenter
$\beta_{\mathcal W_{2}}(\mathcal W_{d-2})$.
\end{proof}

\begin{proof}[Proof of Corollary \ref{cor:1}]
Since regularity and transversality are open conditions, any deformation $\W^{\varepsilon}$ of $\W=\W_{2}\boxtimes \W_{k}$ is of the form $\W^{\varepsilon}=\W^{\varepsilon}_{2}\boxtimes\W^{\varepsilon}_{k}$ with $\W_{k}^{\varepsilon}$ regular and transverse to $\W_{2}^{\varepsilon}$ if $\varepsilon$ is small enough.
By composing by a local diffeomorphism we can assume that
$\Delta(\W^{\varepsilon})=\Delta(\W)$.
Since $\Delta(\W)$ is invariant by $\W_{2}$, it is transverse to $\beta_{\W_{2}}(\W_{k})$ and consequently, it is also transverse to $\beta_{\W_{2}^{\varepsilon}}(\W_{k}^{\varepsilon})$. Since the curvature of $\W^{\varepsilon}$ is holomorphic, Theorem \ref{T:1} implies that $\Delta(\W)$ must be invariant by $\W^{\varepsilon}_{2}$.
\end{proof}

\subsection{Invariant discriminant}
Here we will deal with more degenerate components of the discriminant of
a web. The focus is on irreducible components of the discriminant which are totally
invariant and have minimal multiplicity. Our goal is to show that these components
do not appear in the polar set of the curvature. We start by characterizing the defining
equations of webs having discriminant with these properties.

\medskip

Consider a $\nu$-web $\W_{\nu}$ and let $C$ be an  irreducible component of $\Delta(\W_{\nu})$  which is totally invariant by $\W_{\nu}$.
In suitable coordinates $C\cap U=\{w=0\}$  and $\W_{\nu}$ is defined by
$$dw^{\nu}+w^{m}(a_{\nu-1}(z,w)dw^{\nu-1}dz+\cdots+a_{0}(z,w)dz^{\nu})$$
 for some $m\ge 1$.

\begin{lemma}\label{L:irredutivel}
If $\W_{\nu}$ is as above then $\Delta(\W_{\nu}) \ge (\nu -1) C$, and  equality holds if and  only if $m=1$ and  $a_0(z,0) \neq 0$.  Moreover, in this case
$\W_{\nu}$ is irreducible.
\end{lemma}
\begin{proof}
If  $\W_{\nu}$ is irreducible then  we can consider the Puiseux parametrizations
$$\frac{dw}{dz}=\zeta^{j}c_{0}(z)w^{\frac{r}{\nu}}+\cdots,\qquad c_{0}(z)\not\equiv 0,\quad \zeta=e^{2i\pi/\nu},\quad j=1,\ldots,\nu,$$
 of the defining polynomial of $\W_{\nu}$ in $\C((z))[w,\frac{dw}{dz}]$. Then an equation for
 $\Delta(\W_{\nu})$ is
 $$\prod\limits_{i\neq j}\big((\zeta^{i}c_{0}(z)w^{\frac{r}{\nu}}+\ldots)-(\zeta^{j}c_{0}(z)w^{\frac{r}{\nu}}+\ldots)\big)=w^{r(\nu-1)}\big(c_{0}(z)^{\nu(\nu-1)}\prod\limits_{i\neq j}(\zeta^{i}-\zeta^{j})\big)+\ldots$$
and the multiplicity of $w=0$ is $r(\nu-1)$.
 Moreover, since
$$\prod\limits_{j=1}^{\nu}(\zeta^{j}c_{0}(z)w^{\frac{r}{\nu}}+\ldots)=w^{r}c_{0}(z)^{\nu}+\ldots=w^{m}a_{0}(z,w),$$
we deduce that $ m\le r$. Consequently the multiplicity of $w=0$ is at least $m(\nu-1)\ge \nu-1$.
If $C=\{w=0\}$ has multiplicity $\nu-1$ in $\Delta(\W_{\nu})$ then  $m=r=1$  and $w\not|a_{0}$.

Reciprocally, if $m=1$ and $a_0(z,0) \neq 0$ then $\W_{\nu}$ is irreducible by Eisentein's criterium
and the argument shows that $\Delta(\W_{\nu})= (\nu-1) C$.

\medskip

To conclude the proof of the lemma it suffices to show that when $\W_{\nu}$ is non-irreducible
the inequality $\Delta(\W_{\nu}) > (\nu-1) C$ holds true.
If $\W_{\nu}=\W_{\nu_{1}}\boxtimes\cdots\boxtimes\W_{\nu_{s}}$ is the decomposition of $\W_{\nu}$ in irreducible factors then
$$\Delta(\W_{\nu})=\prod\limits_{i=1}^{s}\Delta(\W_{\nu_{i}})\prod\limits_{i\neq j}\mathrm{tang}(\W_{\nu_{i}},\W_{\nu_{j}}),$$
so that the multiplicity of $w=0$ is at least $\sum\limits_{i=1}^{s}(\nu_{i}-1)+s(s-1)=\nu+s(s-2)$ which is greater than $\nu-1$ if $s\ge 2$.
\end{proof}

We are now ready to establish the regularity of the curvature along totally invariant
irreducible components of the discriminant. Indeed we show more as we allow to
superpose  the irreducible web with a smooth web transverse to it.

\begin{prop}\label{invnu}
Let $\W_{\nu}$ be a $\nu$-web
and let $C\subset\Delta(\W_{\nu})$ be an irreducible component totally invariant by $\W_{\nu}$ and having minimal multiplicity $\nu-1$. Let $\W_{d-\nu}$ be a smooth $(d-\nu)$-web transverse to $C$. Then the curvature of $\W=\W_{\nu}\boxtimes\W_{d-\nu}$ is holomorphic along $C$.
\end{prop}

\begin{proof}
Let $(U,(z,w))$ be a local coordinate system such that $C\cap U=\{w=0\}$,
$$T\W_{\nu}|U=\{dw^{\nu}+w(a_{\nu-1}(z,w)dw^{\nu-1}dz+\cdots+a_{0}(z,w)dz^{\nu})=0\}$$
and
$T\W_{d-\nu}|U=\{\prod\limits_{j=1}^{d-\nu}(dz+\frac{1}{\nu}g_{j}(z,w)dw)=0\}.$
Let $\pi:\bar U\to U$  be the ramified covering given by $(z,w)=\pi(x,y)=(x,y^{\nu})$. The irreducibility of $\W_{\nu}$ implies that its monodromy group is cyclic and consequently $\pi^*\W_{\nu}$ is totally decomposable.
In fact, $\pi^{*}\W_{\nu}$ is given by
$$y^{\nu(\nu-2)}dy^{\nu}+\bar a_{\nu-1}(x,y)y^{(\nu-1)^{2}}dy^{\nu-1}dx+\cdots+\bar a_{1}(x,y)y^{\nu-1}dy\,dx^{\nu-1}+\bar a_{0}(x,y)dx^{\nu}=0.$$
Since $y\not|\bar a_{0}$ we can write the differential $1$-forms defining $\pi^*\W_{\nu}$ as
$$\omega_{i}:=dx+y^{\nu-2}f(x,\zeta^{i}y)\zeta^{-i}dy,\qquad i=1,\ldots,\nu.$$
The differential $1$-forms defining  $\pi^*\W_{d-\nu}$ are
$$\omega_{\nu+j}:=dx+y^{\nu-1}g_{j}(x,y^{\nu})dy,\qquad j=1,\ldots,d-\nu.$$
Recall from \S\ref{SS:curvatura} that $K(\pi^*\W)=\sum\limits_{i<j<k}d\eta_{ijk}$, where
$\eta_{ijk}$ is the unique $1$-form satisfying $d(\delta_{rs}\omega_{t})=\eta_{ijk}\wedge\delta_{rs}\omega_{t}$ for each cyclic permutation $(r,s,t)$ of $(i,j,k)$ and the function $\delta_{rs}$ is defined by $\omega_{r}\wedge\omega_{s}=\delta_{rs}(x,y) dx\wedge dy$.
We denote by $\varphi_{\ell}(x,y)=(x,\zeta^{\ell}y)$, $\ell=1,\ldots,\nu$ the deck transformations of $\pi$.
We have that
$$K(\pi^{*}\W)=\pi^{*}K(\W)=\frac{1}{\nu}\sum_{\ell=1}^{\nu}\varphi_{\ell}^*\pi^{*}K(\W)=\sum_{i<j<k}d\left(\frac{1}{\nu}\sum_{\ell=1}^{\nu}\varphi_{\ell}^*\eta_{ijk}\right).$$
On the one hand  $\frac{1}{\nu}\sum_{\ell=1}^{\nu}\varphi_{\ell}^*(y^{n}dx)=y^{n}dx$ if and only if $n\equiv 0 \mod \nu$ and $\frac{1}{\nu}\sum_{\ell=1}^{\nu}\varphi_{\ell}^*(y^{n}dy)=y^{n}dy$ if and only if $n\equiv -1 \mod \nu$. On the other hand, if $\eta_{ijk}=A_{ijk}(x,y)dx+B_{ijk}(x,y)dy$ then from Lemma~\ref{curv} follows
that the order of the poles of $A_{ijk}$ along $y=0$ is  $\le\nu-1$ and  $B_{ijk}$ is logarithmic along $y=0$ with constant residue. This fact jointly with the previous remark will imply that $d\eta_{ijk}$ is holomorphic along $y=0$ and consequently $K(\W)$ is holomorphic along $C$.
\end{proof}

\section{Global webs and  Legendre transform}\label{S:legendre}

Now we turn our attention to global $k$-webs of degree $d$ on the projective plane. As it was already mentioned in the introduction these are determined by a section $ \omega$ of $Sym^k \Omega^1_{\mathbb P^2}(d+ 2k)$ having isolated zeros and non-zero discriminant. Dually,
a $k$-web of degree $d$ can also be expressed as a section $X$ of $Sym^k T\mathbb P^2( d-k)$ subjected to the very same conditions as above: isolated zeros and non-zero discriminant .

It follows from Euler's  sequence
\[
0 \to \mathcal O_{\mathbb P^2} \longrightarrow \mathcal O_{\mathbb P^2}(1)^{\oplus 3} \longrightarrow  T{\mathbb P^2} \to 0,
\]
that sections $\omega \in H^0 ( \mathbb P^2, \Omega^1_{\mathbb P^2}(d+2))$ and $X \in H^0(\mathbb P^2, T\mathbb P^2(d-k))$ defining the same foliations can be presented in homogeneous coordinates as
\begin{enumerate}[(a)]
\item a homogeneous vector field with coefficients of  degree $d$
$$X=A(x,y,z)\frac{\partial}{\partial x}+B(x,y,z)\frac{\partial}{\partial y}+C(x,y,z)\frac{\partial}{\partial z}, $$ and
\item a homogeneous 1-form with coefficients of degree $d+1$
$$\omega=p(x,y,z)dx+q(x,y,z)dy+r(x,y,z)dz$$
satisfying $\omega(R)=0$, where $R=x\frac{\partial}{\partial x}+y\frac{\partial}{\partial y}+z\frac{\partial}{\partial z}$.
\end{enumerate}
The relation between $X$ and $\omega$ is given by
\begin{equation}
\label{equ1}
\omega=\imath_{R}\imath_{X}\Omega, \qquad\textrm{where}\quad \Omega=dx\wedge dy\wedge dz.
\end{equation}
It is clear that the 1-dimensional distribution $\langle X\rangle$ on $\mathbb C^3$ is not uniquely determined by the foliation, only the 2-dimensional distribution $\ker\omega= \langle X, R\rangle$ is. As $\omega(R)=0$, there exist homogeneous polynomials $A',B',C'$ of degree $d$ such that
\begin{equation}
\label{equ2}
\omega=A'\alpha+B'\beta+C'\gamma,
\end{equation}
where
\begin{equation}\label{abc}
\alpha=ydz-zdy,\quad \beta=zdx-xdz,\quad \gamma=xdy-ydx.
\end{equation}
From (\ref{equ1}) and (\ref{equ2}) it follows that
$$(p,q,r)=(A',B',C')\times (x,y,z)=(A,B,C)\times (x,y,z),$$
so that we can take $(A',B',C')=(A,B,C)+\lambda (x,y,z)$ for any homogeneous polynomial $\lambda$ of degree $d-1$.

\medskip

From Euler's sequence we can deduce the following exact sequence
\[
0 \to Sym^{k-1} ( \mathcal O_{\pp^2} (1)^{\oplus 3} ) \otimes \mathcal O_{\pp^2} \to  Sym^{k} ( \mathcal O_{\pp^2} (1)^{\oplus 3} ) \to Sym^k T\pp^2 \to 0 \, .
\]
It implies that a  $k$-web of degree $d$ on $\pp^{2}$ is determined  by a  bihomogeneous polynomial $P(x,y,z;a,b,c)$ of degree $d$  in the coordinates $(x,y,z)$ and degree  $k$ in the coordinates $(a,b,c)$ respectively. More concretely,
\begin{enumerate}[(a)]
\item $X=P\left(x,y,z;\frac{\partial}{\partial x},\frac{\partial}{\partial y},\frac{\partial}{\partial z}\right)$ determines  a global section of $Sym^{k}T\pp^{2}(d-k)$, and
\item $\omega=P(x,y,z;\alpha,\beta,\gamma)$ determines a global section of $Sym^{k}\Omega^1_{\pp^{2}}(d+2k)$.
\end{enumerate}
Notice that two polynomials $P$ and $P'$ differing by a multiple of $xa+yb+zc$ determine
the same sections.

There exist natural homogeneous coordinates  in the dual projective plane $\check{\pp}^{2}$ which associates to the point
$(a:b:c)\in\check{\pp}^{2}$ the line $\{  ax+by+cz=0\} \subset \pp^2$.
Since
\begin{eqnarray*}
T^*_{(x:y:z)}\pp^{2}&=&\{\omega=a\,dx+b\,dy+c\,dz\in T^*\C^{3}: \omega(R)=0\}\\ &=&\{a\,dx+b\,dy+c\,dz : ax+by+cz=0\}
\end{eqnarray*}
there exists a natural identification of $\pp T^*\pp^{2}$ with the incidence variety
$$\mathcal{I}=\{((x:y:z),(a:b:c)) | ax+by+cz=0\}\subset \pp^{2}\times\check{\pp}^{2}.$$

 \begin{figure}[t!]
  \begin{center}
   \epsfig{file=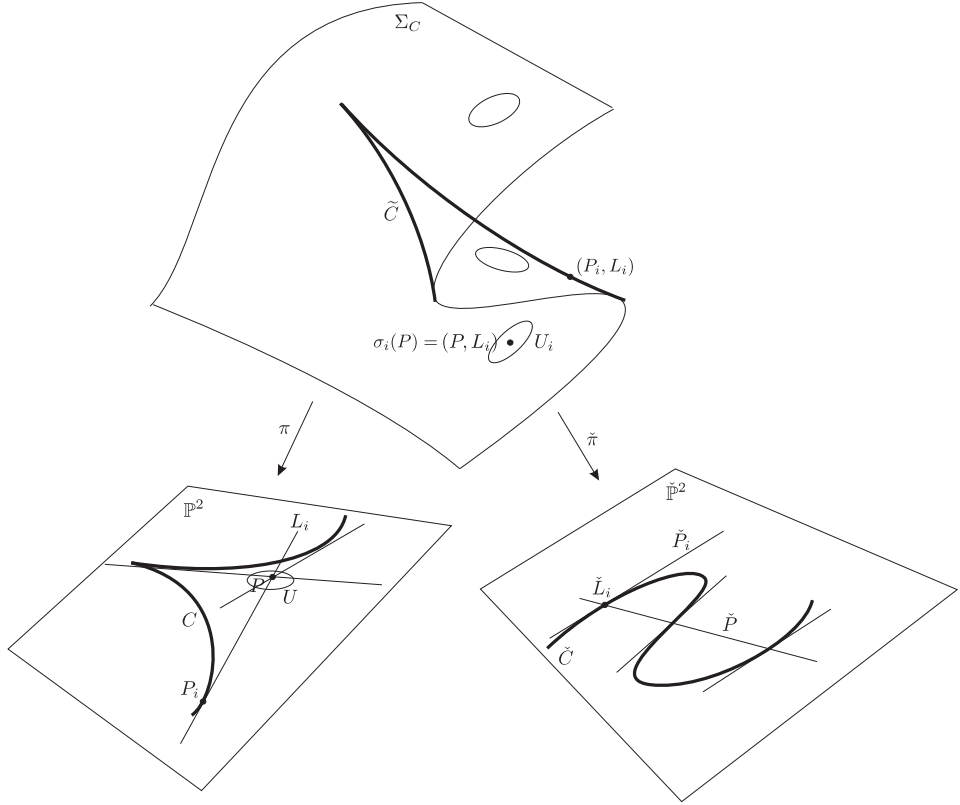,width=\textwidth}
  \end{center}
    \caption{The algebraic web $\leg C$}
  \label{algweb}
 \end{figure}

Let $\mathcal{W}$ be a $k$-web of degree $d$ on $\pp^{2}$ and let $P(x,y,z;a,b,c)$ be a bihomogeneous polynomial defining $\mathcal{W}$.
Then $S_{\mathcal{W}}\subset\pp T^{*}\pp^{2}$, the graph of $\mathcal W$ on $\pp T^* \pp^2$,   is given by
$$S_{\mathcal{W}}=\{((x:y:z),(a:b:c))\in\pp^{2}\times\check{\pp}^{2} | ax+by+cz=0, P(x,y,z;a,b,c)=0\}$$
under the above identification between $\mathcal{I}$ and $\pp T^{*}\pp^{2}$.

Suppose $\mathcal W$ is an irreducible web of degree $d>0$ and consider the restrictions  $\pi$ and $\check{\pi}$ to $S_{\mathcal{W}}$ of the natural projections of $\pp^{2}\times\check{\pp}^{2}$ onto $\pp^{2}$ and $\check{\pp}^{2}$ respectively. These projections $\pi$ and $\check{\pi}$ are rational maps of  degrees $k$ and $d$ respectively. The contact distribution $\mathcal{D}$ on $\pp T^{*}\pp^{2}$ is identified with
$$\mathcal{D}=\ker(a\,dx+b\,dy+c\,dz)=\ker(x\,da+y\,db+z\,dc).$$
The foliation $\F_{\mathcal{W}}$ induced by $\mathcal{D}$ on $S_{\mathcal{W}}$ projects through $\pi$ onto the $k$-web $\mathcal{W}$ and it projects through $\check{\pi}$ onto a $d$-web $\check{\mathcal{W}}$ on $\check{\pp}^{2}$.

\begin{definition}
If $d>0$ and the $k$-web $\W$ of degree $d$ is irreducible then the $d$-web $\check{\mathcal{W}}$ on $\check{\pp}^{2}$ is called the \textbf{Legendre transform} of $\mathcal{W}$ and it  will be denoted by $\leg\mathcal{W}$.
\end{definition}

If $\mathcal W$ is determined by $P\left(x,y,z;\frac{\partial}{\partial x},\frac{\partial}{\partial y},\frac{\partial}{\partial z}\right)$, or respectively by $P(x,y,z;ydz-zdy,zdx-xdz,xdy-ydx)$, then its Legendre transform $\leg \W$ is determined by $P\left(\frac{\partial}{\partial a},\frac{\partial}{\partial b},\frac{\partial}{\partial c};a,b,c\right)$, respectively by $P(bdc-cdb,cda-adc,adb-bda;a,b,c)$.\\

Using these formulae we can proceed to define the Legendre transform for arbitrary $k$-webs of  arbitrary degree $d$. Notice that when $\mathcal W$ decomposes as the product of two webs $\mathcal W_1 \boxtimes \mathcal W_2$ then its Legendre transform will be the product of $\leg(\mathcal W_1)$ with $\leg(\mathcal W_2)$.  But the Legendre transform of an irreducible  $k$-web of degree $0$ is no longer a web on $\pp^2$. It is instead an irreducible curve of degree $k$. Similarly, the Legendre transform of  a reduced curve of degree $d$ is a $d$-web of degree $0$, see Figure 1.

If we consider the space of $k$-webs of degree $d$,  $\mathbb W(k,d) \subset \mathbb P H^0(\pp^2, Sym^k \Omega^1_{\pp^2}(d+2k))$, as
the projectivization of the space of $k$-symmetric $1$-forms with non-zero discriminant
and having singular set with {\bf reduced} divisorial components (instead of only isolated singularities), then
the Legendre transform
defines an involutive isomorphism
\[
\leg : \mathbb W(k,d) \longrightarrow \mathbb W(d,k) \,
\]
when $k,d \ge 0$ and $d+k>0$. 

It is easy to check the following  properties of the Legendre transform:
\begin{enumerate}
[(a)]
\item Let us fix a generic line $\ell$ on $\pp^{2}$. Then $tang(\mathcal{W},\ell)=p_{1}+\ldots+p_{d}$, where $p_{i}\in\pp^{2}$. We can think $\ell$ as a point of $\check{\pp}^{2}$ and the $p_{i}$ as straight lines on $\check{\pp}^{2}$ passing through the point $\ell$. Then $T_{\ell}\leg\mathcal{W}=\bigcup\limits_{i=1}^{d}T_{\ell}p_{i}$.
\item If  $L$ is a leaf of $\mathcal{W}$ distinct from a line then the  union of lines tangent to
$L$ is a leaf of $\leg(\mathcal W)$.
\end{enumerate}

Consider an affine chart $(x,y)$ of $\pp^{2}$ and an affine chart of $\check{\pp}^{2}$ whose coordinates $(p,q)$ correspond to the line $\{y=px+q\}\subset\pp^{2}$. If a web $\W$ is defined by an implicit affine equation $F(x,y;p)=0$ with $p=\frac{dy}{dx}$ then $\leg(\W)$ is defined by the implicit affine equation
\begin{equation}
\label{affine}\check{F}(p,q;x):=F(x,px+q;p)=0,\qquad\textrm{with}\qquad x=-\frac{dq}{dp}.
\end{equation}
In particular, for a foliation  defined by a vector field $A(x,y)\px+B(x,y)\py$  we can take $F(x,y;p)=A(x,y)p-B(x,y)$.\\

We will proceed to describe some of the geometry of the Legendre transform of a foliation
$\F$ on $\pp^2$.  We start by describing the role played by the inflection divisor of $\F$.

\subsection{Inflection divisor for foliations}\label{SS:inf}

Let $\F$ be a degree $d$, $d>0$, foliation of $\mathbb P^2$
 and $X$ any degree $d$ homogeneous vector
field on $\C^3$ inducing $\F$. The inflection divisor of $\F$,
denoted by $I(\F)$, is the divisor defined by the vanishing of the 
determinant
\begin{equation}\label{ext1}
\det \left( \begin{array}{ccc}
  x & y & z \\
  X(x) & X(y)  & X(z) \\
  X^2(x) & X^2(y)  & X^2(z)
\end{array} \right) \, .
\end{equation}

In \cite{extactic}, $I(\F)$ was called the {\em first extactic} curve of
$\F$ and the following properties were proven:
\begin{itemize}
\item[(a)] If the determinant (\ref{ext1}) is identically zero then $\F$ admits a
rational first integral of degree $1$; that is,  if we suppose that
the singular set of $\F$ has codimension $2$ then the degree of
$\F$ is zero;
\item[(b)] On $\mathbb P^2 \setminus Sing({\F})$, $I(\F)$ coincides with the curve described by the
inflection points of the leaves of $\F$;
\item[(c)] If $C$ is an irreducible algebraic invariant curve of $\F$ then
$C \subset I(\F)$ if, and only if, $C$ is an invariant line;
\item[(d)] The degree of $I(\F)$ is exactly $3d$.
\end{itemize}

As a consequence of property  (d) we obtain that the maximum
number of invariant lines for a degree $d$ foliation is $3d$. This
bound is attained, even if we restrict to real foliations and
real lines, as the Hilbert modular foliation of degree $5$
described in Section \ref{S:hilb5} shows.

One can also define the inflection divisor for an arbitrary $k$-web $\mathcal W$ on $\mathbb P^2$.
One has to consider the surface $S_{\mathcal W} \subset \mathbb P T\mathbb P^2$ naturally associated
to $\mathcal W$;  and take the tangency locus $T$  of $S_{\mathcal W}$ with the foliation on $\mathbb P T\mathbb P^{2}$
induced by the lifting of all the lines of $\mathbb P^2$. The inflection divisor of $\mathcal W$ can be
then defined as $\pi_* T$, where $\pi: S_{\W}  \to  \mathbb P^2$ is the natural projection.
Since we will not use it in what follows, we will not provide more details but instead redirect the
interested reader to \cite[Example 2.13]{Maycol}. Here we will just mention that it is a divisor of degree $k^{2}+(2d-1)k+d$.

\medskip

Let $C$ be an irreducible curve contained in the support of the inflection divisor $\mathcal I(\F)$
of $\F$. If $C$ is $\F$-invariant then $C$ is a line and the corresponding point on the dual projective plane is a singular point for $\leg \F$. If instead $C$ is not $\F$ invariant then
the image of $C$ under the Gauss map of $\F$ is a curve $D$ of the dual projective plane
$\check \pp^2$. In general $D$ is not invariant by $\leg \F$. When it is invariant
one has strong implications in the geometry of $\F$ as is stated below.

\begin{prop}\label{D}
If $D$ is $\leg \F$ invariant then we are in one of the two following cases:
\begin{enumerate}
\item The curve $D$ is a line in $\check \pp^2$ and the tangent line of $\F$ at a generic point $p$
of $C$ is the line joining $p$ and the point of $\pp^2$ determined by the line $D$;
\item The curve $D \subset \check \pp^2$  is not a line, its dual curve $\check D \subset \pp^2$ is $\F$-invariant and the tangent line at a generic point of it is
tangent to $\F$ at some point of $C$. Moreover at a neighborhood of a generic point of $D$ the Legendre transform of $\F$ decomposes as the product of foliation tangent to $D$ and a $(d-1)$-web transverse to $D$.
\end{enumerate}
\end{prop}
In particular, for  webs on the projective plane of degree one, on each irreducible component $C$ of the discriminant which is not a straight line, the multiple directions can not be  tangent to $C$ at generic points.

\begin{proof}
If $D$ is a line invariant by $\leg \F$ then the point  $p \in \pp^2$ determined by it must be a singular point of $\F$. Moreover, since $D$ is the image of $C$ under
the Gauss map of $\F$  the generic line  through $p$ must be tangent to $\F$ at some point of $C$.

If $D$ is not a line then tangent lines of $D$ determine the dual curve $\check D \subset \pp^2$. It is clear that $\check D$ must be $\F$-invariant, and
as not every point of $\check D$ is an inflection point, over a generic point of $D$  only one of the tangent lines of $\leg \F$ is tangent to $D$.
\end{proof}

Besides the components of $\Delta(\leg \F)$ determined by the inflection divisor of $\F$ there are also the ones determined
by singularities of $\F$.

\subsection{Singularities versus invariant lines}\label{sing}

If $p \in \pp^2$  is a singularity of a foliation  $\F$ then the line determined by it on $\check \pp^2$ must be invariant by $\leg \F$.
The dual line of a general singularity $p$ will  not be  contained in the discriminant of $\leg \F$. This  will be the case
if and only if the tangency at $p$  between $\F$ and a generic line through $p$ has order at least two or $I(\F)$ contains a non invariant irreducible component of the tangency locus between $\F$ and the pencil of lines through $p$.
This last eventuality does not occur when considering convex foliations.
One can promptly verify that the first eventuality holds if and only if the singularity has zero linear part or if its linear part is a non-zero multiple of the radial vector field. Singularities in the latter situation
will be called {\bf radial singularities}. Note that these, by definition, have non-zero linear part.

Although any two radial singularities are locally analytically equivalent by a classical theorem of Poincar\'e, they may behave distinctly under
the Legendre transform. The point is that a generic line through a radial singularity has tangency with the foliation of multiplicity
at least two but it may be bigger. If we write $X = c_{\nu} R  +  X_{\nu} + h.o.t.$ with $c_{\nu}(0,0) \neq 0$ and $X_{\nu}$ homogenous of degree
$\nu$ and not proportional to $R$ then the generic line through zero has tangency of multiplicity $\nu$ at zero with the foliation determined by $X$.
In this case we will say that the radial singularity has order $\nu-1$. Radial singularities of order one, will be also called {\bf simple radial singularities}.

\begin{prop}\label{rad33}
If $s$ is a radial singularity of order $\nu-1$ of a foliation $\F$ then at a neighborhood of a generic point of the line $\ell$ dual to $ s$
the web $\leg \F$ can be written as the product $\W_1 \boxtimes \W_2$, where $\W_1$ is an irreducible $\nu$-web leaving
$\ell$ invariant and $\W_2$ is a web transverse to $\ell$. Moreover, still at a neighborhood of a generic point of $\ell$,
\[
\Delta(\leg \F) = (\nu-1) \ell + \Delta(\W_2) \, .
\]
\end{prop}
\begin{proof}
By using (\ref{affine}), if $s=(0,0)$ is a radial singularity of order $\nu-1$ of $\F$ then
$$\check{F}(p,q;x)=q+a_{1}(p,q)qx+\cdots+a_{\nu-1}(p,q)qx^{\nu-1}+a_{\nu}(p,q)x^{\nu}+\cdots+a_{d}(p,q)x^{d},$$
with $a_{\nu}(p,0)\not\equiv 0$.
By applying Weierstrass preparation theorem we can write
 $$\check{F}(p,q;x)=U(p,q;x)(x^{\nu}-q(\bar a_{\nu-1}(p,q)x^{\nu-1}+\cdots+\bar a_{0}(p,q)),$$
where $U(p,0;0)\not\equiv 0$.
Hence, $\leg\F=\W_{1}\boxtimes\W_{2}$ near the generic point of  $\ell=\{q=0\}$,
 where $\W_{2}$ is a $(d-\nu)$-web transverse to $\ell$ and $\ell$ is totally
invariant by $\W_{1}$.  Lemma \ref{L:irredutivel} implies that $\W_1$ is irreducible and
$\Delta(\W_1) = (\nu-1) \ell$.
\end{proof}

Let $r_i(\F)$ denote the number of radial singularities of a foliation $\F$ having order $i$.
As the discriminant of a $d$-web of degree $1$ has degree $(d+2)(d-1)$  the proposition
above has the following consequence.

\begin{cor}
If $\F$ is a foliation of degree $d$ then
\[
\sum_i i \cdot r_i( \mathcal F) \le (d+2)(d-1) \, .
\]
\end{cor}

Combining Theorem \ref{T:1} with the previous consideration we obtain a characterization of flat $3$-webs of degree $1$ satisfying some conditions. If $C$ is a  non invariant irreducible  component of the inflection divisor of a degree $3$ foliation $\F$ on $\pp^{2}$ then we consider the curve $C^{\perp}$ consisting of those points $q$ for which there exists  $p\in C$ such that $tang(\F,T_{p}\F)=2p+q$.

\begin{prop}\label{degree3}
Let $\F$ be a degree $3$ foliation on $\pp^{2}$ with reduced inflection divisor $I(\F)$. A necessary condition for $\leg\F$ being flat is that for each non invariant irreducible component $C$ of $I(\F)$ we have that $C^{\perp}$ is invariant by $\F$.
Moreover, if all the singularities of $\F$ have non zero linear part, this condition is also sufficient.
\end{prop}

\begin{proof} First we will prove that for each non invariant irreducible component $C$ of $I(\F)$, $K(\leg\F)$ is holomorphic along $D\subset\Delta(\leg\F)$ if and only if $C^{\perp}$ is invariant by $\F$, where $D$ is the image of $C$ by the Gauss map of $\F$. The reducedness of $I(\F)$ implies that in a neighborhood of $D$ we can decompose $\leg\F=\W_{1}\boxtimes\W_{2}$, where $\W_{2}$ is a $2$-web with discriminant $D$ and
 $\W_{1}$ is a foliation transverse to $\W_{2}$. By Theorem \ref{T:1}, the curvature of $\leg\F$ is holomorphic along $D$ if and only if $D$ is invariant by either $\W_{2}$ or by $\beta_{\W_{2}}(\W_{1})=\W_{1}$. In the first case, we are in the eventuality (1) of Proposition \ref{D} and consequently $C^{\perp}$ is a singular point of $\F$. In the second case, $C^{\perp}$ is contained in the envelope of the family of lines $\{T_{p}\F, p\in C\}$ and consequently it is invariant by $\F$.

Secondly, if all the singularities of $\F$ have non zero linear part then $\Delta(\leg\F)$ only contains the previous considered components $D$ and the dual lines $\ell$ to radial singularities of order $\nu-1\in\{1,2\}$. By Proposition \ref{rad33}, we can decompose  $\leg\F=\W_{\nu}\boxtimes\W_{3-\nu}$ in a neighborhood of $\ell$, with $\ell$ totally $\W_{\nu}$-invariant and $\W_{3-\nu}$ transverse to $\ell$.
From Proposition \ref{invnu} we deduce that $K(\leg\F)$ is holomorphic along $\ell$.
\end{proof}

\section{Duals of convex foliations are flat}

{\bf Convex foliations} are those without inflection points along the leaves which are not straight lines,
i.e. those whose inflection curve is totally invariant (a product of lines).

\subsection{Singularities of convex foliations}

\begin{lemma}\label{L:sing}
Let $\F$ be a convex foliation on $\mathbb P^2$. If the inflection curve
of $\mathcal F$ is reduced then the singularities of $\F$ have non-nilpotent linear
part.
\end{lemma}
\begin{proof}
Let $p$ be a singularity of $\F$ with nilpotent linear part. Throughout we will
assume that $p=(0,0) \in \C^2$. If $X$ is a polynomial vector field inducing $\F$, decompose
it as $X = X_1+ X_2 +\ldots + X_k$ where $X_i$ is a vector with polynomial homogeneous components
of degree  $i$.

The lemma will follow from  a simple analysis of the first non-zero jet at the origin  of
the polynomial
\[
I(X) = \det \left(
  \begin{array}{cc}
    X(x) & X(y) \\
    X^2(x) & X^2(y) \\
  \end{array}
\right)
\]
which defines the inflection curve of $\F$ in $\C^2$. The key observation is that under the hypothesis of
 complete invariance and reduceness of $I(\F)$, the directions determined by the
first non-zero jet of $I(X)$  determine the $\F$-invariant lines through $p$.

Let $X_i$ be the first non-zero jet of $X$. The $\F$-invariant lines  through the
origin  must be invariant by $X_i$.

If  $X_i$ is not proportional to the radial
vector field $R$ then  the $\F$-invariant lines through the origin
are contained in the zero locus of $xX_i(y) - yX_i(x)$. Therefore there
are at most $i + 1$ of them.  If we write down the homogeneous components
of $I(X)$, we promptly realize that
\[
I(X) = X_i(x)X_i^2(y) - X_i(y)X_i^2(x) + h.o.t.
\]
In particular, the algebraic multiplicity of $I(X)$ at $0$ is at least  $3i - 1$. But $3i-1 > i+1$ unless
$i=1$. When this is the case, $xX_1(y) - yX_1(x)$ cuts out the eigenspaces of matrix $DX_1$. The reducedness
of $I(\F)$ excludes the possibility of a nilpotent linear part.

It remains to treat the case where $X_i$ is a multiple of the radial vector field. Let now $X_j$, for some $j>i$,
be the first jet not proportional to $R$. Now the $\F$-invariant lines are in the zero locus of $xX_j(y) - yX_j(x)$.
Thus there are at most $j+1$ of them. On the other hand,  all the  jets of $I(X)$ of order strictly smaller than
$2i + j -1$ are zero. Hence, the reducedness  of $I(\F)$ ensures that the number of  $\F$-invariant
lines through $0$ is at least $2i+j-1$. As before, $2i+j-1 > j+1$ unless $i=1$. This settles the lemma.
\end{proof}

\subsection{Flatness of reduced convex foliations} We are now ready to prove Theorem \ref{T:convex}. We
restate it thinking on  reader's convenience.

\begin{thm}
The dual web of a reduced convex foliation  is flat.
\end{thm}
\begin{proof}
As we have seen in \S  \ref{SS:inf} and \ref{sing}   the discriminant of the web dual to a foliation $\F$ is composed by the dual of its inflection curve
and   lines dual to  some of its singularities. If $\F$ is convex and has degree $d$ then the inflection curve is entirely composed by $3d$ invariant lines and hence its dual is a finite number of points.
Hence $\Delta(\check{\F})$ is a product of lines corresponding to some of the singularities of $\F$.

As explained in \S \ref{sing}, the only singularities of $\F$ that contribute to the discriminant of $\check{\F}$ are those
with zero or radial linear part.  Lemma \ref{L:sing} excludes the case of zero linear part. Let us analyze  a radial singularity $s$ of order $\nu$ of $\F$. The  dual line $\check{s}\subset\check{\mathbb{P}}^{2}$ is contained in the discriminant of $\check{\F}$. Proposition \ref{rad33} implies that  we can decompose locally $\check{\F}=\W_{1}\boxtimes\W_{2}$  near $\check{s}$, where $\W_{1}$ is an irreducible $\nu$-web leaving $\check{s}$ invariant and $\W_{2}$ is a $(d-\nu)$-web transverse to $\check{s}$. We claim that $\W_{2}$ is regular near $\check{s}$, i.e. through a generic point of $\check{s}$ we have $(d-\nu)$ different tangents lines to $\W_{1}$. Indeed, the tangent lines to $\check{\F}$ are the dual of the tangency points of $\F$ with a generic line through $s$. This tangency locus is composed by $s$ itself with multiplicity $\nu$ and other $d-\nu$ points. If two of these points would coincide then we would have an inflection point of $\F$ on a non invariant straight line, which contradicts the convexity of $\F$.
By applying Proposition~\ref{invnu} we deduce that the curvature of $\check{\F}$ is holomorphic along $\check{s}$.
Since $s$ is an arbitrary radial singularity of $\F$, the $2$-form $K(\check{\F})$ is holomorphic on the whole $\check{\mathbb{P}}^{2}$. But $\mathbb P^2$ does not have holomorphic $2$-forms, and therefore the curvature of $\check{\F}$ vanishes identically.
\end{proof}

\section{Examples of convex foliations}\label{S:examples}

In this section we exhibit some examples ensuring that Theorem \ref{T:convex} is not versing about the empty set.
We start by describing an infinite family of convex foliations -- the Fermat foliations --, showing their
birational invariants, and  studying the  algebraization of its
elements. Next we describe three sporadic examples: the Hesse pencil $\mathcal H_4$, the Hilbert modular foliation
$\mathcal H_5$, and $\mathcal H_7$   a foliation of degree $7$ sharing with $\mathcal H_4$ twelve invariant lines.

\subsection{Fermat foliations}
Let $d\ge 2$ be an integer. Consider the foliation $\mathcal G_d$ determined by the rational function
\begin{align*}
g_d : \pp^2 &\dashrightarrow \pp^1 \\
(x: y : z) &\longmapsto \frac{x^{d-1} - y^{d-1}}{z^{d-1} - x^{d-1}} \, .
\end{align*}
The rational function $g_d$ has  three singular fibers $g_d^{-1}(0),g_d^{-1}(\infty), g_d^{-1}(1)$ which are products
of $(d-1)$ lines, and every other fiber is isomorphic to the Fermat curve of degree $d-1$. The foliation $\G_d$
has degree $2d -4$.

Consider now the standard Cremona involution $\varphi : \pp^2 \dashrightarrow \pp^2$, $\varphi(x:y:z) = (x^{-1}:y^{-1}:z^{-1})= ( yz: xz: xy)$. Define $\F_d$ as the pull-back of $\G_d$ by $\varphi$. Thus $\F_d$
is determined  by the rational function
\[
f_d(x:y:z) = \varphi^* h_d(x:y:z) = \frac{(yz)^{d-1} - (xz)^{d-1} } {(xy)^{d-1} - (yz)^{d-1} } = \frac{z^{d-1}(y^{d-1} - x^{d-1})}
{y^{d-1}(x^{d-1} - z^{d-1})} \, .
\]
Each of the fibers $f_d^{-1}(0),f_d^{-1}(\infty), f_d^{-1}(1)$ is a product of $d-1$ reduced lines
and one line with multiplicity $d-1$. The degree of $\F_d$ is
equal to $ 4(d-1) - 2 - 3(d-2) = d$ according to Darboux's formula \cite{Jouanolou}. Since $\F_d$ leaves invariant $3d$ lines it follows that it is a convex foliation.

If $\pi_d: S \to \pp^2$ is the  blow-up at the $(d-1)^2$ base points of $g_d$ then $\hat{g_d} = g_d\circ \pi_d$ is
an isotrivial fibration and the foliation $\pi_d^* \F_d$ is completely transverse to the
smooth fibers of $\hat{g_d}$.

\begin{remark}\rm
If $\F$ is a  foliation on a projective surface $S$ and
$\overline{\F}$ in $\overline{S}$ is any reduced foliation
birationally equivalent to $\F$,  then the
 {\em foliated genus} of $\F$  is
defined  in \cite{LG,LGtese} as
\[ g(\F) := \chi(\oo_{\overline{S}}) + \frac{1}{2}T^* {\overline{\F}} \cdot( T^*{\overline{ \F}}
 \otimes K_{\overline{S}}^*)  \, .\]

If the morphism $\pi_d: S_d \to \mathbb P^2$ corresponds to the
blow-up of the radial singularities of  $\mathcal F_d$ then the
foliation $\overline{\mathcal F_d} = \pi_d^* \mathcal F_d $ is
such that (cf \cite[Example 3.5.1]{LG})
\[
   g(\mathcal F_d) = g(\overline{\mathcal F_d})= \frac{d(d+1)}{2}
   - (d-1)^2 -3  = \frac{ - d^2 + 5d - 8 } {2}\,
\]
and in particular
\[
 \lim_{d \to \infty} g(\F_d) = - \infty \, ,
\]
showing that there is no lower bound for the foliated genus of holomorphic
foliations. This answers a question raised in \cite{LGtese}.
\end{remark}

\subsection{Algebraization of $\check{\F}_d$} The $d$-webs $\check{\F}_d$ are not
only flat but also algebraizable. Indeed they belong to bigger family of algebraizable
webs that we now proceed to describe.


Given $(p,q)$ coprime integers with
$q>0$, we define $\mathcal{F}_{p/q}$ as the $q^2$-web  induced by
$$\omega_{p/q}=  \prod_{n=1}^q\prod_{m=1}^q
(x^{\frac{p}{q}}\alpha+e^{\frac{2i\pi
m}{q}}y^{\frac{p}{q}}\beta+e^{\frac{2i\pi
n}{q}}z^{\frac{p}{q}}\gamma), $$
where $\alpha,\beta,\gamma$ are the homogeneous $1$-forms introduced in (\ref{abc}).
It has degree $d$ where,
\begin{equation}
\label{deg}
d=\left\{\begin{array}{rcl}pq & \textrm{if} & p>0\\ -2pq & \textrm{if} &
p<0\end{array}\right.
\end{equation}
When $q=1$, we recover the $1$-webs (foliations) $\mathcal F_d$.

Given $(p,q)$ coprime integers as above, define the \emph{triangular-symmetric curve} (the terminology is classical, see
for instance \cite{Feld}) of type $(p,q)$ as follows:
$$\mathbb{F}_{p/q}=\prod_{n=1}^q\prod_{m=1}^q (x^{\frac{p}{q}}+e^{\frac{2i\pi
m}{q}}y^{\frac{p}{q}}+e^{\frac{2i\pi n}{q}}z^{\frac{p}{q}}).$$
It follows from  $\gcd(p,q)=1$ that $\mathbb{F}_{p/q}$ is irreducible. More informally we can think of $\mathbb F_{p/q}$ as the curve cut out by the  algebraic function $x^{p/q} + y^{p/q} + z^{p/q}.$

 For $\varepsilon\in\mathbb{Q}^*$,
let us consider the correspondence (multi-valued algebraic map)
$h_{\varepsilon}:\mathbb{P}^{2}\dasharrow\mathbb{P}^{2}$ given by
$$h_{\varepsilon}(x:y:z)=(x^{\varepsilon} : y^{\varepsilon} : z^{\varepsilon}) .$$

\begin{prop}\label{P:fermat}
The following assertions hold true:
\begin{itemize}
\item[(a)] If  $\varepsilon=\frac{1}{q-p}$ then
$h_{\varepsilon}^*\mathcal{F}(p/q)$  is an algebraic web (the
product of $q^{2}$ pencils of lines);
\item[(b)] If $\varepsilon=\frac{p}{p-q}$ then $h_{\varepsilon}^*\mathrm{Leg}(\mathcal{F}(p/q))=\mathrm{Leg}(\mathbb{F}(p/q))$. In
particular the dual of $\mathcal{F}(p/q)$ is an algebraizable $d$-web, where $d$ is given by (\ref{deg}).
\end{itemize}
\end{prop}
\begin{proof} The proof is a blind computation. First notice that
$$h_{\varepsilon}^*\alpha=\varepsilon (yz)^{\varepsilon-1}\alpha,\qquad
h_{\varepsilon}^*\beta=\varepsilon (xz)^{\varepsilon-1}\beta,\qquad
h_{\varepsilon}^*\gamma=\varepsilon (xy)^{\varepsilon-1}\gamma.$$
Then, for $\varepsilon=\frac{1}{q-p}$, one has the following identity:
$$h_{\varepsilon}^*\omega_{p/q}=\varepsilon^{q^{2}}(xyz)^{\frac{pq^{2}}{q-
p}}\prod_{m=1}^{q}\prod_{n=1}^{q}(\alpha+e^{\frac{2i\pi
m}{q}}\beta+e^{\frac{2i\pi n}{q}}\gamma),$$
which implies (a).

\bigskip

Similarly, for $\varepsilon=\frac{p}{p-q}$, one has
\begin{eqnarray*}h_{\varepsilon}^*\mathrm{Leg}(\mathcal{F}(p/q))&=&\varepsilon
^{pq}(xyz)^{\frac{pq^{2}}{p-q}}\prod_{m=1}^{q}\prod_{n=1}^{q}(\alpha^{\frac{p}
{q}}+e^{\frac{2i\pi m}{q}}\beta^{\frac{p}{q}}+e^{\frac{2i\pi
n}{q}}\gamma^{\frac{p}{q}})\\
&=&\varepsilon^{pq}(xyz)^{\frac{pq^{2}}{p-q}}\mathrm{Leg}(\mathbb{F}(p/q)) \, ,
\end{eqnarray*}
which implies (b).
\end{proof}

\subsection{Hessian pencil of cubics}
Every nonsingular cubic in $\pp^2$ is projectively equivalent to one
defined in projective coordinates $(x:y:z)$ by
\[
  F_{\alpha}(x,y,z) = x^3 + y^3 + z^3 - 3 \alpha x y z \, \, , \text{ where
  } \alpha \in \C \text{ and } \alpha^3 \neq 1.
\]
If $H(F)$ denotes the Hessian determinant of $F$ then a direct
computation shows that the solutions of
$F_{\alpha}=H(F_{\alpha})=0$ do not depend on $\alpha \in \C
\setminus \lbrace z | z^3=1\rbrace$. The solutions corresponds to
 the inflection points of $F_{\alpha}$ and there are nine of
them. These nine points lie on $12$ projective lines which are
the four singular cubics described by the parameters $\alpha^3=1$
and $\alpha=\infty$.

Consider the foliation $\mathcal H_4$ induced
by the projective $1$--form $\omega = f dg - g df$ where
$f = x^3+y^3+z^3$ and $g=3 x y z$. It
has degree $4$  and $12$ invariant lines. The
radial singularities of $\omega$ corresponds to the inflection
points of $F_{\alpha}$. Through each of them passes four invariant lines,
and therefore they all have order two.

If $\sigma:S \to \pp^2$ is the blow up of
$\pp^2$ on these radial singularities then $\sigma^* \mathcal H_4$ is a reduced
foliation and its cotangent bundle is $\mathcal O_S(C)$ where $C$ is
the strict transform of one of the cubics of the pencil. It follows that
$\mathrm{kod}(\mathcal H_4) = \nu (\mathcal H_4) = 1$.

The dual of $\mathcal H_4$ is a $4$-web of degree $1$ with discriminant divisor
supported on the union of  the nine lines of $\check \pp^2$ dual to the radial singularities of $\mathcal H_4$.
Since they are all radial singularities of order two, each of these lines
appear  in the discriminant divisor with multiplicity two.

For any $4$-web $\W$ on a complex surface $S$  we can consider the meromorphic map $j: S \dashrightarrow \pp^1$ which sends a point $p$ of $S$ to the $j$-invariant of the four
tangent directions of $\W$ at $p$.  Recall that the $j$-invariant is the unique analytic invariant
 of four unordered points of $\pp^1$. It is for unordered points what the cross-ratio is for ordered points.

The property of the discriminant of $\leg \mathcal H_4$ alluded to above implies that the
$j$-invariant of $\leg \mathcal H_4$  is  identically zero.
Indeed, for any web the $j$-invariant has polar set contained in the discriminant. But, it
is well known, that the $j$-invariant of $3p + q$ with $p\neq q$ is zero. Hence the $j$-invariant
of $\leg \mathcal H_4$ is a  meromorphic function on $\pp^2$ without poles and equal to
zero when restricted to a number of lines. Hence it must be identically zero.

We can apply the proposition below to deduce that $\leg \mathcal H_4$ is a parallelizable, and consequently an algebraizable,  $4$-web.

\begin{prop}
Let $\W$ be a germ of smooth  $4$-web on  $(\C^2,0)$. If the $j$-invariant of $\W$ is
constant then
\[
K(\W) = 4 K (\W')
\]
where $\W'$ is any $3$-subweb of $\W$.
 Moreover if $K(\W) =0$ then $\W$ is parallelizable.
\end{prop}
\begin{proof}
Since $\W$ is a germ of smooth $4$-web it is the product of $k$ distinct foliations $\F_1, \ldots, \F_4$ which we will assume to be defined by $1$-forms $\omega_1, \ldots, \omega_4$.
We can further assume that $\omega_1 + \omega_2 + \omega_3=0$.

The $j$-invariant of $\W$ is constant if and only if a holomorphic multiple of $\omega_4$ can be written as a linear combination of $\omega_1$ and $\omega_2$ with constant coefficients. Thus we
can assume  $\omega_4 =\lambda \omega_1 + \mu \omega_2$ where $\lambda, \mu \in \C$.
This is sufficient to show the existence of a unique holomorphic $1$-form $\eta$ satisfying
\[
d\omega_i = \eta \wedge \omega_i
\]
for every $i = 1, \ldots, 4$. Its differential is the curvature of any $3$-subweb of $\W$. Hence
$K(\W) = 4 K(\W')$ as wanted.

If $K(\W)=0$ then $\eta$ is closed and
so are the $1$-forms $\beta_i= exp\left( \int \eta \right) \omega_i$. The map
\[
(x,y) \mapsto \left( \int \beta_1, \int \beta_2 \right)
\]
conjugates the $4$-web $\W$ to the $4$-web defined by the levels of the linear forms
$ x, y, x+y, \lambda x + \mu y$. Hence $K(\W)=0$ implies $\mathcal W$ parallelizable.
\end{proof}

\begin{remark}
An analogous result holds for $k$-webs ($k > 4$) if one assumes that the $j$-invariant
of any $4$-subweb is constant. The proof is the same.
\end{remark}

\subsection{Hilbert modular foliation}\label{S:hilb5} Our next example
was studied in \cite{LG}. It is
$\mathcal H_5$,  the degree $5$ foliation of $\mathbb P^2$ induced by the following $1$-form  on $\mathbb C^2$:
\[
  (x^2-1)(x^2- (\sqrt{5}-2)^2)(x+\sqrt{5}y) dy - (y^2-1)(y^2- (\sqrt{5}-2)^2)(y+\sqrt{5}x) dx \, .
\]
It  leaves invariant an arrangement of $15$ real lines that can be synthetically
described as follows. Consider the icosahedron embedded in $\mathbb R^3$ with its
center of mass at the origin. Use radial projection to bring it to the unit sphere $S^2$.
On the quotient $\mathbb P^2_{\mathbb R}$ of $S^2$
by the antipodal involution, the $30$ edges of the icosahedron will become $15$ line segments.
The corresponding line arrangement, more precisely an arrangement isomorphic to it,  is left invariant by $\mathcal H_5$.

The foliation $\mathcal H_5$ has $10$ radial singularities of order one, coming from the centers of the 20 faces of the icosahedron,
and $6$ radial singularities of order three, coming from the twelve vertices of the icosahedron.

It has negative Kodaira dimension, and numerical Kodaira dimension equal to one.  See \cite{LG} for a thorough discussion.

Theorem \ref{T:convex} implies that the Legendre transform of $\mathcal H_5$ is flat. A computer-assisted calculation
shows that its linearization polynomial \cite{HenautLinear} has degree four (see remark below), and hence it is not linearizable, and in particular it is not algebraizable.
Indeed we do believe that  $\leg(\mathcal H_5)$ has no abelian relation at all, but so far we do not   have a proof.

\begin{remark}
The linearization polynomial
is not intrinsically attached to a web, one has to choose local coordinates and write it as an implicit differential equation. The claim
about its degree means that in suitable coordinates it has degree four. Anyway this is sufficient to
ensure the non-linearizability of the web.
\end{remark}

\subsection{Degree $7$ foliation invariant by the Hessian
Group} The group of symmetries of the Hessian configuration of
$12$ lines was determined by Jordan as a subgroup of $PSL(2,\C)$
of order $216$. It is generated by projective transformations of
order $3$ which leave one of the $12$ projective lines pointwise
fixed. It also contains nine involutions which fix nine invariant
lines. Together with the twelve lines of the Hessian arrangement
these  nine  lines form an arrangement of $21$ lines, which we will be called
the  {\em extended  Hessian arrangement}, following \cite{Or}.

It is possible to prove that
the degree $7$ foliation $\mathcal H_7$ given in affine coordinates by
\[
            (x^3-1)(x^3+7 y^3+1) x \partial_x + (y^3-1)(y^3+7x^3+1) y \partial_y
\]
is invariant by the Hessian group and leaves invariant the extended
Hessian arrangement of lines.

It is tangent to a pencil of curves of degree $72$. Except for three special elements,
the generic member of the pencil has genus $55$. The special elements are:
\begin{enumerate}
\item a completely decomposable fiber, with support equal to the extended Hesse arrangement.
The $12$ irreducible components appearing in the original Hesse arrangement have multiplicity $3$,
while the remaining $9$ appear with multiplicity $4$;
\item a fiber of multiplicity three, with support equal to an irreducible curve of degree $24$ and
genus $19$;
\item a fiber of multiplicity two, with support equal to an irreducible curve of degree $36$ and
genus $28$.
\end{enumerate}
These claims have been verified with the help of a computer.

The foliation $\mathcal H_7$ carries $21$ radial singularities: $12$ with multiplicity  $3$, and $9$ with multiplicity $4$.
It has Kodaira dimension and numerical Kodaira dimension equal to two.

Theorem \ref{T:convex} implies that the Legendre transform of $\mathcal H_7$ is also flat. Its
linearization polynomial \cite{HenautLinear} has degree six and, as the Legendre transform of $\mathcal H_5$, it is not algebraizable.
In contrast we do know that  $\leg(\mathcal H_7)$ has at least three linearly independent  abelian relations coming from
holomorphic $1$-forms on a ramified covering  of $\pp^1$  (Klein's quartic) but we
do not know what is the exact rank of it.

\section{Deformations of radial singularities and convex foliations}

The remaining of the paper is devoted to the proof of Theorem \ref{T:3}. The starting point is the following result which
guarantees the persistence of simple radial singularities when we deform a reduced convex foliation in such a way that its
Legendre transform is still flat.

\begin{thm}\label{T:rad}
Let $\F^{\varepsilon}$ be a small analytic deformation of
foliations of  degree $d\ge 2$ on $\mathbb{P}^{2}$.
Suppose that $s^{0}\in\pp^{2}$ is a simple radial singularity of $\F^{0}$ and assume
that the tangency locus between $\F^0$ and the pencil of lines through $s^{0}$ does not contain any non-invariant irreducible component of $I(\F^{0})$, for instance this is the case if $\F^{0}$ is convex. If the dual webs $\leg\F^{\varepsilon}$ are flat then there exists an analytic curve $\varepsilon\mapsto s^{\varepsilon}$
such that $s^{\varepsilon}$ is a simple radial singularity of $\F^{\varepsilon}$.
\end{thm}

\begin{proof}
Since the tangency locus between $\F^0$ and the pencil of lines through $s^{0}$ does not contain any non-invariant irreducible component of $I(\F^{0})$ it follows that
the dual web $\check{\mathcal F}^0$ decomposes at a neighborhood of a general point of $C^{0}$ ( the line dual to $s^{0}$ )   as $\W_{2}^{0}\boxtimes\W_{d-2}^{0}$, where $\W_{d-2}^{0}$ is a smooth web transverse to $C^{0}$.
Taking an affine chart $(x,y)$ in $\mathbb{P}^{2}$ such that $s^{0}=(0,0)$ and the corresponding affine chart $(p,q)$ in $\check{\mathbb{P}}^{2}$, we have that $\check{\F}^{0}$ is given by
$F^{0}(p,q;x)=\sum\limits_{i=0}^{d}a_{i}^{0}(p,q)x^{i}=0$ with $q|a_{i}^{0}$ for $i=0,1$ and $q^{2}\not\!|\,a_{2}^{0}$, so that $C^{0}=\{q=0\}$ is a reduced invariant component of $\Delta(\check{\F}^{0})$.  By continuity, there exists an irreducible  component $C^{\varepsilon}\subset\Delta(\check{\F}^{\varepsilon})$ deforming $C^{0}$.
By Corollary~\ref{cor:1}, locally
$\check{\F}^{\varepsilon}=\W_{2}^{\varepsilon}\boxtimes\W_{d-2}^{\varepsilon}$ and $C^{\varepsilon}=\Delta(\W_{2}^{\varepsilon})$ is invariant by $\W_{2}^{\varepsilon}$.
Thus,  eventuality (2) of Proposition \ref{D} is not possible and $C^{\varepsilon}$ must necessarily be a straight line in $\check{\pp}^2$, dual of some point $s^{\varepsilon}\in\mathbb{P}^{2}$. Since $\check{s}^{\varepsilon}$ is invariant by
$\check{\F}^{\varepsilon}$ we deduce that $s^{\varepsilon}$ must be a singularity of $\F^{\varepsilon}$. Taking into account the discussion of section \ref{sing} we obtain that $s^{\varepsilon}$ is a radial singularity of $\F^{\varepsilon}$.
Taking affine charts $(p_{\varepsilon},q_{\varepsilon})$ in $\check{\mathbb{P}}^{2}$ such that $C^{\varepsilon}=\{q_{\varepsilon}=0\}$ we can present $\check{\F}^{\varepsilon}$
by an equation
$$F^{\varepsilon}(p_{\varepsilon},q_{\varepsilon};x_{\varepsilon})=\sum\limits_{i=0}^{k}a_{i}^{\varepsilon}(p_{\varepsilon},q_{\varepsilon})x_{\varepsilon}^{i}=0,$$
where $q_{\varepsilon}|a_{i}^{\varepsilon}$ for $i=0,1$. By continuity, $q_{\varepsilon}\not\!|\,a_{2}^{\varepsilon}$ if $\varepsilon$ is small enough.
Therefore, $\F^{\varepsilon}$ is given by a vector field
$$c_{0}^{\varepsilon}(x_{\varepsilon}\partial_{x_{\varepsilon}}+y_{\varepsilon}\partial_{y_\varepsilon})+X_{2}^{\varepsilon}+
\cdots$$ where $X_{2}^{\varepsilon}$ is an homogeneous vector field of degree $2$ in the variables $(x_{\varepsilon},y_{\varepsilon})$ not collinear with $x_{\varepsilon}\partial_{x_{\varepsilon}}+y_{\varepsilon}\partial_{y_\varepsilon}$ because this is so when $\varepsilon=0$. Notice that for small $\varepsilon$ we have that $c_{0}^{\varepsilon}\neq 0$ if $c_{0}^{0}\neq 0$.
\end{proof}

\begin{cor}\label{C:radiais}
Let $\F^{\varepsilon}$ be an analytic deformation of the foliation $\F^{0}:=\F_{d}$, $d \ge 3$, such that $\check\F^{\varepsilon}$ is flat for all $\varepsilon\approx 0$. Then $\F^{\varepsilon}$ has at least $(d-1)^2$ simple radial
singularities.
\end{cor}
\begin{proof}
Let $p_1, \ldots, p_{(d-1)^2}$ be the  singularities of $\F_d$ defined by $x^{d-1}-y^{d-1} = x^{d-1} - z^{d-1} = x^{d-1} - y^{d-1} =0$. Since through each of them there are only three $\F_d$-invariant lines, the convexity of
$\F_d$ implies that each of these singularities is radial of order one.
Theorem \ref{T:rad} implies the existence of  $(d-1)^2$ simple radial singularities
for $\F^{\varepsilon}$.
\end{proof}

\section{Rigid flat webs I: the rational case}
In this section
we will study the deformations of $\check{\mathcal F}_3$, the next section will be devoted to the
 deformations of $\check{\mathcal F}_d$ for $d \ge 4$.

\begin{thm}\label{T:rigrational}
The closure of the orbit by $\mathrm{Aut}(\mathbb{P}^{2})$ of the dual web of the foliation $\F_{3}$ is an irreducible component of the space of flat $3$-webs of degree $1$.
\end{thm}

The foliation $\F_{3}:(x^{3}-x)\frac{\partial}{\partial x}+(y^{3}-y)\frac{\partial}{\partial y}$ has $6$ hyperbolic singularities, $4$ radial singularities of order $1$ and $3$ radial singularities of order $2$.
The dual of each radial singularity of order $1$ of $\F_{3}$ is an invariant reduced component of the discriminant of the dual web $\check{\F}_{3}$. The main ingredient in the proof of Theorem \ref{T:rigrational} is the stability of radial singularities of order $1$ given by Theorem \ref{T:rad}.

\subsection{\bf Flat deformations of $\F_3$ are Riccati}
Let $\mathcal G$ be a foliation tangent to a pencil of rational curves on $\pp^2$, and let $\rho:S \to \mathbb P^2$
be a morphism for which $\rho^* \mathcal G$ is a fibration  $\pi:S \to \mathbb P^1$. We will say that a foliation $\F$ is Riccati with respect to $\mathcal G$ if  $\rho^*\F$ has no tangencies with the generic fiber of $\pi$.
In these circumstances  there exists a Zariski open set $U \subset \pp^1$ such that every fiber over a point of $U$
is transverse to $\rho^* \mathcal F$. Moreover, once a base point  $b \in U$ is chosen, there is a natural
representation
\[
\varphi : \pi_1(U,b) \longrightarrow \mathrm{Aut}( \pi^{-1} (b)) \simeq 	\mathrm{Aut}(\pp^1)
\]
called the monodromy representation, obtained by lifting paths in $U$ along the leaves of $\rho^* \mathcal F$.

\begin{lemma} The foliation $\F_3$ is Riccati with respect to $\F_{-1}$. The open set $U$ can be taken equal
to the complement of $3$ points in $\pp^1$, and  the monodromy representation is a morphism from the free group
with two generators onto a subgroup of $\mathrm{Aut}(\pp^1)$ isomorphic to $\mathbb Z_2 \times \mathbb Z_2$.
\end{lemma}
\begin{proof}
Notice that $\F_{-1}$ is the pencil of conics through the points $[\pm 1 : \pm 1 : 1]$.
Let $\rho:S \to \pp^2$ be the blow-up of these four points and $E_1, \ldots, E_4$ be the exceptional
divisors. The foliation $\rho^* \F_3$ has cotangent bundle isomorphic to
\[
T^* \rho^*  \mathcal F_3 = \rho^* T^* \F_3 \otimes \mathcal O_S(-E_1 - \ldots - E_4) = \rho^* \mathcal O_{\pp^2}(2) \otimes \mathcal O_S(-E_1 - \ldots - E_4)
\]
and the strict transform $C$ of a conic through the four $\rho(E_i)$ is defined by a section of
the same line bundle. On the one hand,
\[
\left(T^* \rho^*  \mathcal F_3\right) ^2=T^* \rho^*  \mathcal F_3  \cdot C = C^2 =  0 \, .
\]
On the other hand,  $ T^* \rho^*  \mathcal F_3  \cdot C = C^2 - \mathrm{tang}( \rho^*  \mathcal F_3,   C)$
for any curve $C$ not invariant by $\rho^*  \mathcal F_3$ according
to \cite[Proposition 3, Chapter 3]{Br1}. It follows that $\mathrm{tang}( \rho^*  \mathcal F_3,   C)=0$. Hence $\F_3$ is Riccati
with respect to $\F_{-1}$.

The fibers of the fibration $\pi:S \to \pp^1$  determined by $\F_{-1}$ which are  not completely transverse to $\rho^*\F_3$ are precisely the three singular fibers of $\pi$. On $\pp^2$ they correspond to the six invariant
lines of $\F_{-1}$ which are also invariant by $\F_3$. The other three  $\F_3$-invariant lines intersect
a curve of the pencil of conics in two distinct points away from the base locus. They correspond to
orbits of order two of the monodromy representation. The generic leaf of $\F_3$ is a quartic with smooth
points  at the base locus of the pencil of conics. Hence its strict transform intersects $C$ at four distinct points.
It follows that the image of the monodromy representation has order four. Putting all together we deduce  that the image of the monodromy
representation is a subgroup of $\mathrm{Aut}(\pp^1)$ conjugated to the one generated by $x\mapsto -x$ and $x \mapsto x^{-1}$.
\end{proof}

\begin{lemma}\label{L:Riccati}
Let $\F^{\varepsilon}$ be an analytic deformation of the foliation $\F^{0}:=\F_{3}$ such that $\check\F^{\varepsilon}$ is flat for all $\varepsilon\approx 0$. Then there exists a family $g^{\varepsilon}$ of automorphisms
 of $\mathbb P^2$ such that $(g^{\varepsilon})^*\mathcal F^{\varepsilon}$ is a Riccati foliation with respect to
 $\F_{-1}$. Moreover the tangency between $(g^{\varepsilon})^*\mathcal F^{\varepsilon}$ and $\F_{-1}$ is equal to the
 six $\F_2$-invariant lines.
\end{lemma}
\begin{proof}
Let $\F^{\varepsilon}$ be an analytic deformation of the foliation $\F^{0}:=\F_{3}$ such that $\check\F^{\varepsilon}$ is flat for all $\varepsilon\approx 0$. Since the $4$ radial singularities of order $1$ of $\F^{0}$ are in general position, and they are stable by deformation by Theorem \ref{T:rad}, we can assume that the four points $(\pm 1,\pm 1)$ are also radial singularities of $\F^{\varepsilon}$ of order $1$.
As in the previous lemma one can show that $\F^{\varepsilon}$  is Riccati with respect to $\F_{-1}$. Moreover, as a line through a radial singularity $p$ of $\F^{\varepsilon}$
 has local tangency of order at least two, the six lines joining the four points $(\pm 1,\pm 1)$
must be invariant by  $\F^{\varepsilon}$. As they are also invariant by $\mathcal F_{-1}$
they must contained in $\mathrm{tang}(\F^{\varepsilon},\F_{-1})$. Since the tangency divisor
of foliations of degree $d_1$ and $d_2$ has degree $d_1+d_2+1$, the lemma follows.
\end{proof}

\subsection{\bf Proof of Theorem \ref{T:rigrational}}
On the one hand $\F^{\varepsilon}$  is a transversely affine foliation because its dual is flat, on the other hand   Lemma \ref{L:Riccati} implies $\F^{\varepsilon}$ is a Riccati foliation.  According to a result of
Liouville, see for instance \cite{Frank}, a Riccati foliation is transversely affine if and only if there
exists an invariant algebraic curve generically transverse to the fibration. Consequently the  monodromy of $\mathcal F^{\varepsilon}$ must have a periodic orbit.

Suppose  one of the generators of the monodromy of $\F^{\varepsilon}$, say
the one deforming $ x \mapsto -x$, has non constant conjugacy class. We can assume
that it takes the form $x \mapsto -\lambda(\varepsilon) x$ for some germ of non constant holomorphic function $\lambda$. Since the only points of $\pp^1$ with finite orbit
under   $x \mapsto -\lambda(\varepsilon) x$ for generic $\varepsilon$ are $0$ and $+\infty$,
the other generator of the monodromy, in this same coordinate for $\pp^1$, must be of the form
$x \mapsto \mu(\varepsilon) x^{-1}$ for a suitable germ of holomorphic function $\mu$. But
these are clearly  conjugated to $x \mapsto x^{-1}$. Therefore the conjugacy class of
the local monodromy around at least two of the three invariant fibers do not vary. Consequently
the analytical type of the singularities of $\F^{ \varepsilon}$ on these fibers are the same as
the ones for $\F_3$. In particular $\F^{\varepsilon}$ has at least two extra radial singularities, and
the line $\ell$ joining them is invariant by it.

Consider now the inflection curve of $\F^{\varepsilon}$. If it contains an irreducible component $D$ which is not $\F^{\varepsilon}$-invariant then, by  applying Proposition \ref{degree3}, either (a) the tangents of $\F^{\varepsilon}$ along $D$ intersect at a singular point
of $\F^{\varepsilon}$; or (b) the tangents of $\F^{\varepsilon}$ along $D$ are also tangent to a $\F$-invariant curve $D^{\perp}$ of degree
at least two.

If we are in case (a) then there exists a singularity $p$ of $\F^{\varepsilon}$ for which the tangency divisor  of the pencil of lines $\mathcal L_{p}$ through $p$ and $\F$ vanishes along $D$ with multiplicity two.  Since this holds for every $\varepsilon\neq 0$ small, it follows that  $\F_3$ has a singularity $p$ for which the tangency between $\F_3$ and $\mathcal L_p$, the radial foliation singular at $p$,  is
a non-reduced divisor $T$ of degree $4$. As through the radial singularities of $\F_3$ passes three distinct invariant lines,   $p$ cannot be radial. If $p$ is a reduced singularity then $T$  would have support at two invariant lines through $p$ and a non-invariant line. As the tangency locus between $\mathcal F$ and $\mathcal L_p$  must contain all the singularities of $\F$ and each invariant line contains exactly $4$ singularities,   the non-invariant line would contain $6$ singularities of $\F$. This leads to a contradiction, a non-invariant line contains at most $3$ singularities, which shows that situation (a) is not possible.

If we are in case (b) then $D^{\perp}$ is a $\F^{\varepsilon}$-invariant curve
distinct from $\ell$ and  the six $\F_{-1}$-invariant lines. Hence its strict transform  is invariant by $\rho^* \F$
and it is generically transverse to the fibers of the fibration defined by $\F_{-1}$.
It follows that the monodromy group of $\rho^*\F$ has two distinct periodic orbits. Consequently $\lambda(\varepsilon)$ must be constant and the analytical type of the singularities of $\F^{\varepsilon}$ do not vary with $\varepsilon$. At this point we can see that $\F^{\varepsilon}$ has $7$ radial singularities and share with $\F_3$ the same $9$ invariant lines.
In particular the tangency of $\F_3$ and $\F^{\varepsilon}$ has degree at least $9$.
To conclude one has just to observe that the tangency divisor between two distinct degree $3$ foliations has degree $7$.
\qed

\section{Rigid Flat Webs II: Elliptic and Hyperbolic cases}\label{S:rigid2}

We will now prove Theorem \ref{T:3} when $d\ge 4$. Indeed, according to Corollary \ref{C:radiais}, more will be done as
we will characterize deformations of $\F_d$ for which the $(d-1)^2$ radial singularities
of order one persist.
For $d=4$ there is a pencil of foliations with this property. This pencil
has been studied before by Lins Neto in \cite{LN_Exemplo}. Our result below shows that
there are no other non-trivial deformations up to homographies.

\begin{thm}\label{d4}
Let $\F_{4}^{\varepsilon}$, $\varepsilon\in(\C,0)$ be a deformation of the Fermat foliation $\F_{4}$. If each $\F_{4}^{\varepsilon}$ has $9$ radial singularities then there exists a family of homographies $h_{\varepsilon}\in\mathrm{PGL}(3,\C)$ and an analytic germ  $f:(\C,0)\to(\C,0)$ such that $h_{\varepsilon}^*\F_{4}^{\varepsilon}$ is defined by the family of vector fields $Z_{0}+f(\varepsilon)Z_{1}$, where $Z_{0}=(x^{3}-1)x\partial_{x}+(y^{3}-1)y\partial_{y}$ defines $\F_{4}$ and $Z_{1}=(x^{3}-1)y^{2}\partial_{x}+(y^{3}-1)x^{2}\partial_{y}$ defines the Fermat foliation $\F_{-2}$.
\end{thm}

Using a Maple script by Ripoll we verified  that every element of this pencil of degree $4$ foliations give rises to a
flat $4$-web of degree one and that the generic element of the pencil is not algebraizable, unlike $\leg(\F_{-2})$ and $\leg(\F_4)$.
Indeed there are only $8$ algebraizable elements in the pencil, four of them isomorphic to $\F_4$ and the other four isomorphic to $\F_{-2}$.
 It would be nice to give a geometric proof of these facts.

\medskip

When $d\ge 5$, the foliation $\F_d$ does not admit non-trivial deformations  preserving the $(d-1)^2$ radial singularities of order one.

\begin{thm}\label{d5}
Let $\F_{d}^{\varepsilon}$, $\varepsilon\in(\C,0)$, be a deformation of the Fermat foliation $\F_{d}$, $d\ge 5$. If each $\F_{d}^{\varepsilon}$ has $(d-1)^{2}$ radial singularities then the deformation is analytically trivial, i.e. there exists a family of homographies $h_{\varepsilon}\in\mathrm{PGL}(3,\C)$ such that $\F_{d}^{\varepsilon}=h_{\varepsilon}^*\F_{d}$.
\end{thm}

\begin{lemma}\label{ab}
If a polynomial $F(x,y)$ of degree $\le d$ belongs to the ideal generated by $x^{d-1}-1$ and $y^{d-1}-1$ then there exist affine polynomials $\alpha(x,y)$ and $\beta(x,y)$ such that $F(x,y)=\alpha(x,y)(x^{d-1}-1)+\beta(x,y)(y^{d-1}-1)$.
\end{lemma}

\begin{proof}
Let $\alpha',\beta'$ be polynomials such that $F=\alpha'(x^{d-1}-1)+\beta'(y^{d-1}-1)$. Define $m=\max(\deg\alpha',\deg\beta')$ and consider the homogeneous part $\alpha_{m}'$ of $\alpha'$ of degree $m$. If $n\ge d-1$ then $y^{n}=(y^{d-1}-1)q_{n}(y)+r_{n}(y)$, $\deg(r_{n})\le d-2$. Thus, there exists a suitable polynomial $\kappa$ such that the homogeneous part $\alpha_{m}$ of $\alpha=\alpha'+\kappa(y^{d-1}-1)$ of degree $m$ is of the form
$$\alpha_{m}=\alpha_{m,0}x^{m}+\alpha_{m,1}x^{m-1}y+\cdots+\alpha_{m,d-2}x^{m-d+2}y^{d-2}.$$
Change $\alpha'$ by $\alpha$ and $\beta'$ by $\beta=\beta'-\kappa (x^{d-1}-1)$.
If $m>1$ then $\alpha_{m}x^{d-1}+\beta_{m}y^{d-1}=0$. Hence $x^{m+1}|\beta_{m}y^{d-1}$ and consequently $\beta_{m}=\alpha_{m}=0$.
\end{proof}

\begin{proof}[Proof of Theorem \ref{d5}]
Since the first part of the proof also applies to the case $d=4$, we  will not restrict to the case $d\ge 5$ unless it is strictly necessary.

Write $X_{\varepsilon}=X_{0}+\varepsilon^{k}X_{1}+\varepsilon^{k+1}X_{2}+\ldots$, a vector field defining $\F_{d}^{\varepsilon}$, where
$$X_{0}=(x^{d}-x)\partial_{x}+(y^{d}-y)\partial_{y},\quad X_{1}=(f+xh)\partial_{x}+(g+yh)\partial_{y},$$ $f,g$ are polynomials of degree $\le d$ and $h$ is an homogeneous polynomial of degree $d$. Since all the singularities of $\F_{d}$ are nondegenerate they are stable. After composing by a family of homographies of $\mathbb{P}^{2}$ we can normalize the deformation $\F_{d}^{\varepsilon}$ so that $\{(1:0:0),(0:1:0),(0:0:1),(1:1:1)\}\subset \mathrm{Sing}(\F^{\varepsilon}_{d})$. This implies that $f(0,0)=g(0,0)=f(1,1)=g(1,1)=0$ and $h(x,y)=xy\hbar(x,y)$ for some homogeneous polynomial $\hbar$ of degree $d-2$.

By assumption, there are $(d-1)^2$ germs of holomorphic maps $$p_{ij}: (\mathbb C, 0) \to \mathbb P^2,\qquad i,j=1,\ldots,d-1,$$
such that $p_{ij}(0)=(\zeta^{i},\zeta^{j})$ and $p_{ij}(\varepsilon)$ is a radial singularity of order one for $\mathcal F^{\varepsilon}_d$, where $\zeta$ is a primitive $d-1$ root of the unity. In fact, the only explicit property of $X_{0}$ that we will use in the sequel is that
\begin{itemize}
\item[($\star$)]
$p_{ij}(0)$ are radial singularities for $X_{0}$.
\end{itemize}

Write $p_{ij}(\varepsilon)=p_{ij}(0)+\varepsilon^{\ell} q_{ij}(\varepsilon)$. A straightforward computation shows that if $q_{ij}(0)\neq 0$ then $\ell\ge k$. In fact, we can take $\ell=k$ with $q_{ij}(0)=-DX_{0}(p_{ij}(0))^{-1}(X_{1}(p_{ij}(0)))$.
Since the matrices $DX_{\varepsilon}(p_{ij}(\varepsilon))$ are diagonal (in fact they are multiple of the identity) we obtain that
$\partial_{y}f+x\partial_{y}h$ and $\partial_{x}g+y\partial_{x}h$ vanish at $p_{ij}(0)$. Hence, there exist polynomials $\alpha,\beta,\gamma,\delta$ such that
\begin{equation}\label{fgh}\left\{
\begin{array}{rcl}
\partial_{y}f+x^{2}(\hbar+y\partial_{y}\hbar)&=&\alpha(x^{d-1}-1)+\beta(y^{d-1}-1)\\
\partial_{x}g+y^{2}(\hbar+x\partial_{x}\hbar)&=&\gamma(x^{d-1}-1)+\delta(y^{d-1}-1)
\end{array}
\right.
\end{equation}

By Lemma \ref{ab} we can assume that  $\alpha,\beta,\gamma$ and $\delta$ are affine.
By equating the homogeneous parts of degree $d$ in (\ref{fgh}) we obtain that $x^{2}(\hbar+y\partial_{y}\hbar)=\alpha_{1}x^{d-1}+\beta_{1}y^{d-1}$. Hence $\beta_{1}=0$ and $\partial_{y}^{2}\hbar=0$. Analogously, $\partial_{x}^{2}\hbar=0$ and
consequently we have
\begin{equation}\label{hbar}
\hbar=\left\{\begin{array}{rcl}0 & \textrm{if} & d\ge 5\\
\lambda xy & \textrm{if} & d =  4,
\end{array}\right.
\end{equation}
for some $\lambda\in\C$.

At this point we will assume that $d\ge 5$. The equality $h=0$ means that the line $z=0$ is invariant by the foliation defined by $X_{1}$.
Notice that the initial foliation $\F^0_d$ is invariant under the following automorphism of $\mathbb P^2$: $(x:y:z) \mapsto (y:z:x)$. Since the lines $x=0$, $y=0$, and $z=0$ are permuted by this automorphisms
 we deduce that the lines $x=0$ and $y=0$ are also  invariant by $X_{1}$. Therefore  $x|f$ and $y|g$. Since $h=0$ and $\deg(\partial_{y}f)\le d-1$, by applying Lemma \ref{ab} we deduce that $\alpha,\beta,\gamma$ and $\delta$ are constant. Therefore,
$$f(x,y)=\alpha(x^{d-1}y-y)+\beta(y^{d}/d-y)+\bar f_{1}(x).$$
Since $x|f$ we have that $\alpha=\beta=0$ and $f(x,y)=\bar f_{1}(x)=xf_{1}(x)$. Analogously, $g(x,y)=yg_{1}(y)$. This means that through the points $(1:0:0)$ and $(0:1:0)$ pass $d+1$ lines invariant  by $X_{1}$. By symmetry, the same property must be true for the point $(0:0:1)$, so the tangency locus $xg(y)-yf(x)=xy(g_{1}(y)-f_{1}(x))=0$ between $X_{1}$ and the radial vector field $x\partial_{x}+y\partial_{y}$ is homogeneous. Therefore $f_{1}(x)=\alpha x^{d-1}+\alpha'$ and $g_{1}(y)=\beta y^{d-1}+\beta'$. Using that $f_{1}(1)=g_{1}(1)=0$ we deduce that $\alpha=-\alpha'=-\beta'=\beta$, so that
$f(x)=\alpha(x^{d}-x)$ and $g(y)=\alpha(y^{d}-y)$, i.e.
$X_{1}=\alpha X_{0}$. Thus, $X_{\varepsilon}=X_{0}+\varepsilon^{k}X_{1}+\cdots=(1+\alpha\varepsilon^{k})X_{0}+\varepsilon^{k+1}X_{2}+\cdots$. We conclude by an inductive argument on $k$.
\end{proof}

\begin{proof}[Proof of Theorem \ref{d4}]
We will use the same notations as in the proof of Theorem \ref{d5}. By (\ref{hbar}), $h(x,y)=\lambda x^{2}y^{2}$ and consequently the line $z=0$ is not invariant by $X_{1}$, but the tangency locus of $X_{1}$ with $z=0$ is $2(1:0:0)+2(0:1:0)$. Interchanging the coordinates $x,y,z$ we also deduce that the tangency locus of $X_{1}$ with $x=0$ (resp. $y=0$) is $2(0:1:0)+2(0:0:1)$ (resp. $2(1:0:0)+2(0:0:1)$). This implies that $f(0,y)$ (resp. $g(x,0)$) is a constant multiple of $y^{2}$ (resp. $x^{2}$).

From (\ref{fgh}) and Lemma \ref{ab} we deduce that $\alpha=2\lambda y+\alpha_{0}$ and $\beta=\beta_{0}$ with $\alpha_{0},\beta_{0}\in\C$. Therefore, $f(x,y)=\alpha_{0}x^{3}y+\frac{\beta_{0}}{4}y^{4}-\lambda y^{2}+\bar f_{1}(x)$. Since $f(0,y)$ is a constant multiple of $y^{2}$ we obtain that $\beta_{0}=\alpha_{0}+\beta_{0}=\bar f_{1}(0)=0$ and consequently, $f(x,y)=xf_{1}(x)-\lambda y^{2}$. Analogously, $g(x,y)=yg_{1}(y)-\lambda x^{2}$.
Since $f_{1}(1)=g_{1}(1)=0$, the points  $p_{3j}(\varepsilon)=(1,\zeta^{j})\mod \varepsilon^{k+1}$ (resp. $p_{i3}(\varepsilon)=(\zeta^{i},1)\mod \varepsilon^{k+1}$) belong to the line $x=1$ (resp. $y=1$) through $(1:0:0)$ (resp. $(0:1:0)$). By symmetry, the points $p_{ii}(\varepsilon)=(\zeta^{i},\zeta^{i})\mod \varepsilon^{k+1}$ belong to the line $y=x$ through $(0:0:1)$. Therefore $f_{1}(x)=\alpha(x^{3}-1)$ and $g_{1}(y)=\beta(y^{3}-1)$. Finally, by imposing that $(1,1)$ is a radial singularity we obtain that $\alpha=\beta$ and
hence $$f+xh=\alpha x(x^{3}-1)+\lambda y^{2}(x^{3}-1),\qquad g+yh=\alpha y(y^{3}-1)+\lambda x^{2}(y^{3}-1).$$
Thus, $X_{1}=\alpha Z_{0}+\lambda Z_{1}$ and
$X_{\varepsilon}=(1+\alpha\varepsilon^{k})Z_{0}+\lambda\varepsilon^{k}Z_{1}+\cdots$. Thanks to ($\star$) we can iterate this procedure taking as $X_{0}$ the vector field $Z_{0}+\lambda\varepsilon^{k}Z_{1}$, obtaining that $X_{\varepsilon}$ is parallel to $Z_{0}+f(\varepsilon)Z_{1}$ for some analytic map $f:(\C,0)\to(\C,0)$.
\end{proof}

\section{Questions}

In this final section we highlight some of the questions that naturally emerged in our investigation.
The first question  concerns the classification of reduced convex foliations. It is a curious fact that
all the examples are invariant by complex reflection groups  and the inflection divisor is supported
on the arrangement of the corresponding reflection lines. We believe that the examples presented in Section
\ref{S:examples}  encompass all the reduced convex foliations. As we are not bold enough to pose this
as a conjecture, we instead propose the following problem.

\begin{problem}
Are there any other reduced convex foliations ?
\end{problem}

\medskip

Our second question appeared already in Section \ref{S:examples} and it can be succinctly
stated as follows.

\begin{problem}
Compute the ranks of the webs $\leg(\mathcal H_5)$ and $\leg(\mathcal H_7)$.
\end{problem}

The interest is not just on the answer but on the methods used to obtain them.
It is our believe, already conjectured in \cite{MPP}, that the existence of abelian
relations for webs implies that the foliations involved have Liouvillian first integrals.
If this is true then the rank $\leg(\mathcal H_5)$ would be zero as it is not a transverselly
affine foliation, see \cite{Hilb}. On the other hand we have no idea how to
determine the rank of $\leg(\mathcal H_7)$. It is not even excluded the possibility of being
an exceptional $7$-web.

\medskip

The next problem seems to be very wild in nature, and a complete answer is probably out of reach already
for pretty small values of $k$ and $d$. Nevertheless,  due to the paucity  of examples of flat webs and exceptional
webs in the literature, the task seems to be worth pursuing.

\begin{problem}
Determine (some of) the irreducible components of the space of flat $k$-webs of degree $d$ on $\mathbb P^2$, for small $k$ and $d$.
\end{problem}

Every $k$-web of degree $0$, being algebraic, is automatically flat. Therefore $\mathbb{FW}(k,0)=\mathbb{W}(k,0)$ for any $k\ge 3$.
Already for $(k,d) = (3,1)$, the first non-trivial case, the task of determining the irreducible components of $\mathbb{FW}(k,d)$ seems to be far from trivial.
Even very particular instances of the above problem seems to have interest. For instance, one can ask if
 Bol's $5$-web, seen as a $5$-web of degree $2$ on $\mathbb P^2$, admits non-trivial flat deformations among the $5$-webs of degree $2$.

\medskip

Another particular instance of the problem above concerns the flat webs $\leg(\mathcal F_{p/q})$ introduced in Section \ref{S:examples}.

\begin{problem}
Determine the flat deformations of the webs $\leg(\mathcal F_{p/q})$ for arbitrary relatively prime integers $p$ and $q$.
\end{problem}

It is an interesting problem already when $q=1$ and $p<0$. In this case we are dealing with the Legendre transforms of
foliations of degree $2p$ defined by the pencils of Fermat curves $\{ \lambda (x^{p+1} - y^{p+1}) + \mu (y^{p+1} - z^{p+1}) =0\}$.
We know that there are deformations of these foliations keeping the $(p+1)^2$ simple radial singularities, see \cite[Example 3.1]{PS}.
As we have seen in Section \ref{S:rigid2}, $\leg(\F_{-2})$ has non trivial flat deformations.
It is possible that something similar holds true for other negative values of $p$.

\end{document}